\documentclass[reqno,11pt,twoside,english]{amsart}
\usepackage[T1]{fontenc}
\usepackage[latin9]{inputenc}
\usepackage{amsthm}
\usepackage{amssymb}
\usepackage{amsmath,amsfonts,epsfig}

\makeatletter

\numberwithin{equation}{section}

\usepackage{color}
\usepackage[final,linkcolor = blue,citecolor = blue,colorlinks=true]{hyperref}

\usepackage{type1cm,ae}
\usepackage{graphicx}
\usepackage{xspace}
\usepackage{setspace}
\usepackage[english]{babel}
\usepackage{tikz}
\usetikzlibrary{arrows,patterns}

 \newcommand{\rr}{{\mathbb R}}
 \newcommand{\R}{{\mathbb R}}
  \newcommand {\F}{{\mathbb F}}

 \newcommand\nn{{\mathbb N}}

 \newcommand\loc{{\mathop\mathrm{\,loc\,}}}

 \newcommand\diam{{\mathop\mathrm{\,diam\,}}}

\newcommand{\mH}{{\mathcal H}}
\newcommand{\bx}{{\partial \Omega}}
\newcommand{\real}{{\mathbb R}}
\def\vint{\mathop{\mathchoice%
		{\setbox0\hbox{$\displaystyle\intop$}\kern 0.22\wd0%
			\vcenter{\hrule width 0.6\wd0}\kern -0.82\wd0}%
		{\setbox0\hbox{$\textstyle\intop$}\kern 0.2\wd0%
			\vcenter{\hrule width 0.6\wd0}\kern -0.8\wd0}%
		{\setbox0\hbox{$\scriptstyle\intop$}\kern 0.2\wd0%
			\vcenter{\hrule width 0.6\wd0}\kern -0.8\wd0}%
		{\setbox0\hbox{$\scriptscriptstyle\intop$}\kern 0.2\wd0%
			\vcenter{\hrule width 0.6\wd0}\kern -0.8\wd0}}%
	\mathopen{}\int}

\numberwithin{equation}{section}
\newtheorem{theorem}{Theorem}[section]

\newtheorem{corollary}[theorem]{Corollary}

\newtheorem{lemma}[theorem]{Lemma}

\newtheorem{proposition}[theorem]{Proposition}

\theoremstyle{definition}
\begingroup
\newtheorem{definition}[theorem]{Definition}
\newtheorem{remark}[theorem]{Remark}

\endgroup

\def\XXint#1#2#3{{\setbox0=\hbox{$#1{#2#3}{\int}$}
		\vcenter{\hbox{$#2#3$}}\kern-.5\wd0}}

\makeatother

\usepackage{babel}
\begin{document}

\title[The Dirichlet Problem for Orlicz-Sobolev mappings between metric spaces]{The Dirichlet Problem for Orlicz-Sobolev mappings between metric space}

\author[W.-J. Qi]{Wen-Juan Qi}

\address[Wen-Juan Qi]{Research Center for Mathematics and Interdisciplinary Sciences, Shandong University 266237,  Qingdao, P. R. China}
\email{wenjuan.qi@mail.sdu.edu.cn}

\thanks{W.-J. Qi was supported by the Young Scientist Program of the Ministry of Science and Technology of China (No.~2021YFA1002200), the National Natural Science Foundation of China (No.~12101362) and the Natural Science Foundation of Shandong Province (No.~ZR2021QA003).}

\begin{abstract}
In this paper, we solve the Dirichlet problem for Orlicz-Sobolev maps between singular metric spaces that extends the corresponding result of Guo et al. [arXiv 2021]. As an intermediate step, we develop a version of Rellich-Kondrachov compactness theorem for Orlicz-Sobolev mappings between metric spaces that extends a previous result of Guo and Wenger [Comm. Anal. Geom. 2020]. Another crucial ingredient is an Orlicz-Sobolev extension of the trace theory for metric valued Sobolev maps developed by Korevaar and Schoen [Comm. Anal. Geom. 1993]. 
\end{abstract}

\maketitle

{\small
\smallskip
\newline
}
\bigskip

\arraycolsep=1pt

\section{Introduction}

The classical Dirichlet problem associated to the harmonic mapping system in an Euclidean domain $\Omega\subset \R^n$ asks for a continuous map $u\colon \Omega\to \R^m$ so that 
\begin{equation*}
	\begin{cases}
		\Delta u&=0 \quad \text{ in } \Omega,\\
		\quad u&=f                 \quad \text{ on } \bx.
	\end{cases}
\end{equation*}
An equivalent way to formulate the Dirichlet problem is to consider energy miniming mappings via the Euler-Lagrange equations. To be more precise, one considers minimizers of the Dirichlet energy 
$$E^2(u):=\int_{\Omega}|\nabla u|^2 dx.$$
When we move from Euclidean spaces to Riemannian manifolds, we also have a natural definition of Dirichlet energy. More precisely, given two Riemannian manifolds $(M,g)$ and $(N,h)$, the Dirichlet energy functional is defined as
$$E^2(u):=\int_M |\nabla u|^2 d\mu,$$
where $|\nabla u|$ is the Riemannian length of the gradient of $u$ and $\mu$ is Riemannian volume induced by $g$ on $M$.  

One of the classical methods to solve the Dirichlet problem in the smooth setting is the direct method from the calculus of variations. To apply it, one needs a \emph{Rellich-Kondrachov compactness theorem} for manifold valued Sobolev maps. Roughly speaking, the Rellich-Kondrachov compactness theorem for Sobolev maps says that if $\{u_k\}_{k\in \mathbb{N}}$ is sequence of Sobolev maps with uniformly bounded Sobolev norm, then up to a subsequence, $u_k$ converges in $L^2$ (indeed even in $L^p$ for all $p<2^*:=\frac{2n}{n-2}$) to some limiting Sobolev map $u$ of the same class. For the theory of Sobolev type spaces and the associated Rellich-Kondrachov compactness theorem in Euclidean spaces, see \cite{Adams-book}. 

Now, consider a mapping $u\colon X\to Y$, where $X=(X,d_X,\mu)$ is a metric measure space and $Y=(Y,d_Y)$ a metric space. Unlike the smooth Riemannian case, there is no natural Dirichlet energy functional associated to a sufficiently regular map. Indeed, there are several well-known (and generally different) energy functionals existing in the literature: the Korevaar-Schoen energy functional, the Hajlasz energy functional \cite{H96}, the upper gradient energy functional \cite{hk98,s00} and so on; see \cite{hkst12} for more energy functionals and the associated Sobolev spaces of metric valued maps. For the solvability of Dirichlet's problem (associated to various different energy functionals) in the setting of singular metric spaces, see \cite{s01,GW20,G21,GHWX21} and the references therein.

In this article, we shall focus on the upper gradient energy functional introduced by Heinonen and Koskela \cite{hk98} and solve the associated Dirichlet problem. Throughout this paper, $X=(X,d_X,\mu)$ is assumed to be a complete metric measure space, $Y=(Y,d_Y)$ a complete metric space, $\Omega\subset X$ a bounded domain and $\mH$ a $\sigma$-finite Borel regular measure on $\bx$. For notational simplicity, we sometimes drop the subscripts $X,Y$ from the distances $d_X, d_Y$ and simply write $d$. Given expressions $a$ and $b$, we write $a\lesssim b$ if there is a constant $C>0$ such that $a\leq Cb$.

Fix an $N$-function $\Phi$, let $N^{1,\Phi}(\Omega,Y)$ be the Orlicz-Sobolev spaces based on upper gradients and $M^{1,\Phi}(\Omega,Y)$ the Orlicz-Sobolev spaces based on Hajlasz gradients (see Section \ref{sec:preliminaries} below for precise definition).  The basic theory of Orcliz-Sobolev spaces on metric measure spaces based on upper gradients has been studied in details in the monograph \cite{H04} and then the theory was greatly extended in \cite{T09,H10,HK22}. 


We first recall the definition of trace for metric valued maps, introduced in \cite{GHWX21}.

\begin{definition}\label{trace-metric}
	Let $u\colon \Omega\rightarrow Y$ be a $\mu$-measurable map. Fix a point $x\in \bx$. If for some point $Tu(x)\in Y$, it holds
	\begin{equation}\label{eq:trace metric definition}
		\lim_{r\rightarrow 0^+}\vint_{B(x, r)\cap\Omega}d_Y(u, Tu(x))\, d\mu=0, 
	\end{equation}
	then we say that the trace $Tu(x)$ of $u$ at $x\in \bx$ exists. Also, we say that $u$ has a trace $Tu$ on  $\bx$ if $Tu(x)$ exists for $\mH$-almost every $x\in \bx$.
\end{definition}

For our purpose, it is convenient to separate a class of admissible domains so that the Dirichlet problem is solvable, which plays a similar role as bounded Lipschitz domains in a smooth Riemannian manifold. The following definition was greatly motivated by \cite[Definition 1.2]{GHWX21}.

\begin{definition}\label{def:admissible domain}
	We say that a bounded domain $\Omega\subset X$ is weakly $(\Phi,\theta)$-admissible for some $N$-function $\Phi\in\Delta'\cap\nabla_2$ and $\theta>0$, if 
	\begin{itemize}
		\item $\mu$ is a doubling measure on $\Omega$;
		\item  $\mH$ is upper codimension-$\theta$ regular on $\partial\Omega$;
		\item $\Omega$ supports a local $\Phi$-Poincar\'e inequality with  $\sum_{k\geq 0}\frac{2^{-k}}{\Phi^{-1}\left(2^{-k\theta}\right)}<\infty$;
		\item $N^{1,\Phi}(\Omega, Y) = M^{1,\Phi}(\Omega, Y)$ with comparable norms for all Orlicz-Sobolev maps.
	\end{itemize}
	We say that $\Omega\subset X$ is $(\Phi,\theta)$-admissible if in addition $\Omega$ supports a global $\Phi$-Poincar\'e inequality, that is, for $u\in N^{1,\Phi}(\Omega)$ with $Tu=0$ $\mH$-almost everywhere on $\bx$, it holds
	\begin{equation}\label{eq:global poincare inequality}
		\|u\|_{L^{\Phi}(\Omega)}\leq C(\Omega)\|g_u\|_{L^\Phi(\Omega)}.
	\end{equation}
\end{definition}

For the next concept, we refer to monograph \cite{GW20} for the notion of a non-principal ultrafilter $\omega$ on $\mathbb{N}$ and the definition of ultra-limit $\lim_{\omega}a_m$ of a bounded sequence $\{a_m\}$ of real numbers. Let $(Y,d)$ be a metric space and $\omega$ a non-principal ultrafilter on $\mathbb{N}$. Denote by $Y_\omega$ the set of equivalent classes $[(y_m)]$ with the sequence $\{y_m\}$ in $Y$ satifying $\sup_md(y_1,y_m)<\infty$, where sequences $\{y_m\}$ and $\{y_m'\}$ are indentified if $\lim_{\omega}d(y_m,y_m')=0$. The metric space obtained by equipping $Y_\omega$ with the distance $d_{\omega}([(y_m)],[(y_m')])=\lim_{\omega}d(y_m,y_m')$ is called the ultra-completion or ultra-product of $Y$ with respect to $\omega$. It is clear that $Y$ isometrically embeds into $Y_\omega$ via the map $\iota\colon Y\to Y_\omega$, which assigns to $x$ the equivalent class $[(x)]$ of the constant sequence $\{x\}$. The following definition was introduced in \cite{GW20}.
\begin{definition}\label{def:1-complemented}
	A metric space $Y$ is said to be 1-complemented in some ultra-completion of $Y$ if there exists a non-principal ultrafilter $\omega$ on $\mathbb{N}$ for which there is a 1-Lipschitz retraction from $Y_\omega$ to $Y$.
\end{definition}
The class of metric spaces that are 1-complemented in some ultra-completion includes all proper metric spaces, all dual Banach spaces, some non-dual Banach spaces such as $L^1$, all Hadamard spaces and injective metric spaces; see \cite[Proposition 2.1]{GW20}.

 For each $u\in N^{1,\Phi}(\Omega,Y)$, let $E^\Phi(u)$ be the $\Phi$-weak upper gradient energy functional of $u$ (see Section \ref{sec:metic valued via upper gradients} below for precise definition). Our first main result establishes a general existence result for the Dirichlet problem associated for Orlicz-Sobolev maps between singular metric spaces.  

\begin{theorem}\label{thm:main existence}
	Suppose $\Omega\subset X$ is a $(\Phi,\theta)$-admissible domain and $Y$ is a metric space that is 1-complemented in some ultra-completion of $Y$. Then for each $\phi\in N^{1,\Phi}(\Omega,Y)$, there exists a mapping $u\in N^{1,\Phi}(\Omega,Y)$ with $Tu=T\phi$ such that 
	\begin{equation*}
		E^\Phi(u)=\inf\left\{E^\Phi(v): v\in N^{1,\Phi}(\Omega,Y)\ \text{ and }\ Tv=T\phi\right\}.
	\end{equation*}	
\end{theorem}

Theorem \ref{thm:main existence} extends the correpsonding result of \cite{GHWX21} from the Sobolev class to more general Orlicz-Sobolev class.


One important ingredient in the proofs is the following version of  Rellich-Kondrachov compactness theorem for Orlicz-Sobolev maps, for which we state as a separate theorem below. It extends a recent result of Guo and Wenger \cite[Theorem 3.1]{GW20}. 

\begin{theorem}\label{thm:main theorem}
	Suppose $\Omega\subset X$ is a weakly $(\Phi,\theta)$-admissible domain for some N-function $\Phi\in\Delta'\cap\nabla_2$. For every $m\in\nn$, let $(Y_m,d_m)$ be a complete metric space, $K_m\subset Y_m$ compact and $\{u_m\} \subset N^{1,\Phi}(\Omega, Y_m)$. Suppose that $(K_m, d_m)$ is a uniformly compact sequence and 
	\begin{equation}\label{eq:bdd-energy-distance-seq-um}
		\sup_{m\in\nn}\left[  \|d_m(y_m,u_m)\|_{L^\Phi(\Omega)}+E^\Phi(u_m)\right]<\infty
	\end{equation}
	for some and thus every $y_m\in K_m$. Then, after possibly passing to a subsequence, there exist a complete metric space $Z$, isometric embeddings $\varphi_m\colon Y_m \hookrightarrow Z$, a compact subset $K\subset Z$ and $v\in N^{1,\Phi}(\Omega, Z)$ such that $\varphi_m(K_m)\subset K$ for all $m\in \nn$ and $\varphi_m\circ u_m$ converges to $v$ in $L^\Phi(\Omega, Z)$.
\end{theorem}

Recall that a sequence of compact metric spaces $(B_m, d_m)$ is called \emph{uniformly compact} if $\sup_{m}\diam B_m<\infty$ and if for every $\varepsilon>0$ there exists $N\in\nn$ such that every $B_m$ can be covered by at most $N$ balls of radius $\varepsilon$. 

Taking $\Phi(x)=|x|^p$, $p>1$ in Theorem \ref{thm:main theorem}, we obtain the corresponding Rellich-Kondrachov compactness theorem for Sobolev spaces.

\begin{corollary}\label{coro:compactness in Sobolev space}
	Suppose $\Omega\subset X$ is a weakly $(p,\theta)$-admissible domain for some $p>1$. For every $m\in\nn$, let $(Y_m,d_m)$ be a complete metric space, $K_m\subset Y_m$ compact and $\{u_m\} \subset N^{1,p}(\Omega, Y_m)$. Suppose that $(K_m, d_m)$ is a uniformly compact sequence and 
	\begin{equation*}
		\sup_{m\in\nn}\left[\int_\Omega d_m^p(y_m, u_m)\, d\mu + E^p(u_m)\right]<\infty
	\end{equation*}
	for some and thus every $y_m\in K_m$. Then, after possibly passing to a subsequence, there exist a complete metric space $Z$, isometric embeddings $\varphi_m\colon Y_m \hookrightarrow Z$, a compact subset $K\subset Z$ and $v\in N^{1,p}(\Omega, Z)$ such that $\varphi_m(K_m)\subset K$ for all $m\in \nn$ and $\varphi_m\circ u_m$ converges to $v$ in $L^p(\Omega, Z)$. 
\end{corollary}

In the setting of singular metric spaces, Rellich-Kondrachov compactness theorem for Sobolev functions has already been considered by Hajlasz-Koskela \cite[Theorem 4]{HaK98}. When $\Omega\subset \R^n$ is a bounded Lipschitz domain, Corollary \ref{coro:compactness in Sobolev space} reduces to the recent theorem of Guo and Wenger \cite[Theorem 3.1]{GW20}. It plays a key role in their solution of the Dirichlet problem for Sobolev mappings with values in certain locally noncompact (infinite dimensional) metric spaces. 

We would like to emphasize that in Theorem \ref{thm:main theorem} (or Corollary \ref{coro:compactness in Sobolev space}), different with the classical Rellich-Kondrachov compactness theorem, we did not claim that a subsequence of $\{u_m\}$ converges, but a subsequence of the modified sequence $\{\varphi_m\circ u_m\}$ converges to some $v$ in $L^p$. However, when $X$ and $Y_m$ are geometrically nice, it is possible to obtain a subsequence of $\{u_m\}$ that do converge; see \cite[Section 3]{GZ21} for such type of results. Indeed, Corollary \ref{coro:compactness in Sobolev space} plays a fundamental role in the compactness results for energy minimizing harmonic mappings between Alexandrov spaces considered in \cite{GZ21} and the Dirichlet problem for $p$-harmonic mappings between singular metric spaces considered in \cite{GHWX21}.

In the formulation of Theorem \ref{thm:main existence}, we need the fact that the trace operator $T$ is well-defined on $N^{1,\Phi}(\Omega,Y)$. When $\Omega\subset X$ is a bounded Lipschitz domain in a smooth Riemannian manifold, $Y$ is a complete metric space and $\Phi(x)=|x|^p,\,1<p<\infty$, this fact was established by Korevaar-Schoen in \cite[Section 12]{ks93}. Very recently, this was extended to Sobolev mappings between metric spaces in \cite{GHWX21}. We give a further extension of this fact to Orlicz-Sobolev maps between metric spaces.

\begin{theorem}\label{thm:trace bdd-1}
	Suppose $\Omega$ is a weakly $(\Phi,\theta)$-admissible domain and $Y$ is a complete metric space embedded isometrically into some Banach space. Then the trace operator 
	$$T\colon N^{1, \Phi}(\Omega, Y, d\mu)\rightarrow L^\Phi(\bx, Y, d\mH)$$ 
	is bounded and linear.
\end{theorem}

Another crucial fact that we shall need in the proof of Theorem \ref{thm:main existence} is the following convergence result for traces of Orlicz-Sobolev maps with uniformly bounded energy.  When $\Omega\subset X$ is a bounded Lipschitz domain in a smooth Riemannian manifold, $Y$ is a complete metric space and $\Phi(x)=|x|^p,\,1<p<\infty$, this fact was established by Korevaar-Schoen in \cite[Theorem 1.12.2]{ks93}, and when $\Phi(x)=|x|^p,\, 1<p<\infty$ it was established recently in \cite[Theorem 1.6]{GHWX21}.  

\begin{theorem}\label{thm:trace convergence1}
	Suppose $\Omega\subset X$ is a weakly $(\Phi,\theta)$-admissible domain and $Y$ is complete. Let $\{u_i\}\subset N^{1,\Phi}(\Omega,Y)$ be a sequence with uniformly bounded energy, that is, 
    $$\sup_{i\in \mathbb{N}}E^{\Phi}(u_i)<\infty.$$
	If $u_i$ converges to some $u\in N^{1,\Phi}(\Omega,Y)$ in $L^\Phi(\Omega,Y)$, then $Tu_i\to Tu$ in $L^\Phi(\bx,Y)$.
\end{theorem}

Many of the arguments in our proofs are similar to the one used in \cite{GW20,GHWX21} and so we do not claim any essential new ideas or techniques produced in this paper. However, as one could expect, since a general Orlicz function (different from the typical $|x|^p$) involves, the estimates often becomes much more delicate. For the convenience of the readers, we have included as many details as possible in all the following proofs.


The paper is organized as follows. In Section \ref{sec:preliminaries}, we introduce the basic notation and collect some auxiliary results. In Section \ref{sec:proof}, we prove the Rellich-Kondrachov compactness theorem both in Orlicz-Hajlasz-Sobolev spaces and in Orlicz-Sobolev spaces based on upper gradients. In Section \ref{sec:trace}, we give an extension of the trace theory of Korevaar-Schoen and prove our trace theorem. In Section \ref{sec:dirichlet problem}, we solve the Dirichlet problem, i.e. Theorem \ref{thm:main existence}.

\section{Preliminaries and auxiliary  results}\label{sec:preliminaries}

Throughout this paper, $X=(X, d_X,\mu)$ is assumed to be a complete metric measure space, $(Y, d_Y)$ a complete   metric space and $\Omega\subset X$ a bounded domain. We say that the measure $\mu$ is a \emph{doubling measure} on $\Omega$ if there exists a constant $C_d\geq 1$ such that 
$$0<\mu(B(x, 2r)\cap\Omega) \leq C_d\,\mu(B(x, r)\cap \Omega)<\infty$$
for all $x\in \bar\Omega$ and $r>0$, where  $B(x, r):=\{y\in X: d(y, x)<r\}$ denotes the open ball centered at $x$ with radius $r$.

Given a set $G\subset \bar\Omega$ endowed with a $\sigma$-finite Borel regular measure $\mH$, we  say that $\mH$ is {\it upper codimension-$\theta$ regular on $G$} for some $\theta>0$ if there exists a constant $C_G$ such that
\begin{equation}\label{upper}
	\mH(B(x, r)\cap G) \leq C_G \frac{\mu(B(x, r)\cap \Omega)}{r^\theta} 
\end{equation}
for  all $x\in G$ and $r>0$.

\subsection{Orlicz spaces}
    In this section, we shall recall the definition of Orlicz spaces. The presentation is rather standard and can be found for instance in \cite{RR91}.
    
	Let $\Phi\colon\rr\to \bar \rr^+$ be a convex function satisfying the conditions: $\Phi(-x)=\Phi(x)$, $\Phi(0)=0$, and $\lim\limits_{x\to\infty}\Phi(x)=+\infty$. With each such function $\Phi$, one can associate another convex function $\Psi\colon\rr\to \bar \rr^+$ having similar properties, which is defined by$$\Psi(y) := \sup\{x|y|-\Phi(x):x\geq 0\}, \quad y\in\rr.$$
Then $\Phi$ is called a \emph{young function}, and $\Psi$ the \emph{complementary function} to $\Phi$. It follows from the definition that $\Psi(0)=0$, $\Psi(-y)=\Psi(y)$ and that $\Psi$ is a convex increasing function satisfying $\lim\limits_{y\to\infty}\Psi(y)=+\infty$. 

Let $\Phi_1, \Phi_2$ be two Young functions. Then $\Phi_1$ is said to be \emph{essentially stronger} than $\Phi_2$, denoted $\Phi_2\prec\prec\Phi_1$, if$$\Phi_2(x)\leq\Phi_1(ax)$$holds for each $a>0$ and for all $x\geq 0$.

Young functions can be classified based on their growth rate. A Young function $\Phi$ is said to be \emph{doubling} or satisfies the \emph{$\Delta_2$-condition}
 if there is a constant $C_2>0$ such that $$\Phi(2x) \leq C_2\Phi(x)$$for each $x\geq 0$. The $\Delta_2$-condition tells that for large $x$ the growth of a Young function $\Phi$ is dominated by the function $C|x|^p$ with some $p>1$ and a constant $C>0$; see \cite[the proof of Corollary 2.3.5]{RR91}. In particular, the
$\Delta_2$-condition excludes functions with exponential growth.

A condition in the opposite direction is the $\nabla_2$-condition. A Young function $\Phi$ satisfies the \emph{$\nabla_2$- condition} if$$\Phi(x) \leq \frac{1}{2C}\Phi(Cx)$$for some fixed constant $C>1$ and for each $x\geq 0$. 

A continuous Young function $\Phi\colon\rr\to\rr^+$ is said to be an \emph{N-function} if $\Phi$ satisfies $\Phi(x)=0$ if and only if $x=0$, 
$\lim\limits_{x\to 0}\frac{\Phi(x)}{x}=0$ and $\lim\limits_{x\to \infty}\frac{\Phi(x)}{x}=\infty$.

An N-function $\Phi$ satifies the \emph{$\Delta'$-condition} if there is a constant $C'>0$ such that $$\Phi(xy)\leq C'\Phi(x)\Phi(y)$$ for all $x,\,y\geq 0$. This is stronger condition than the doubling condition, see \cite[the proof of Lemma 2.3.8]{RR91}. Moreover, $\Delta'\cap\nabla_2\neq \emptyset$, for example, consider the function $\Phi(x)=|x|^\alpha(|\log|x||+1),\,\alpha>1$.

	Let $\Phi$ be a Young function and $(\Omega,\Sigma,\mu)$ be an arbitrary measure space, where $\Sigma$ is a $\sigma$-algebra. The {\it Orlicz space} $L^{\Phi}(\Omega)$ is defined to be $$L^{\Phi}(\Omega)=\left\{u : \Omega\to\rr : u\,\, \text{measurable}, \int_\Omega \Phi(\alpha |u|)\,d\mu<\infty \,\text{for some $\alpha>0$}\right\}.$$
	
	As in the theory of $L^p$-spaces, the elements in $L^\Phi(\Omega)$ are actually equivalence classes consisting of functions that differ only on a set of $\mu$-measure zero. The Orilcz space $L^\Phi(\Omega)$ is a Banach space when equipped with the norm
	$$\Arrowvert u \Arrowvert_{L^{\Phi}(\Omega)}=\inf\left\{k>0: \int_\Omega\Phi\left(\frac{u}{k}\right)\,d\mu\leq 1\right\}.$$

For further analysis, we give an extension of H\"older's inequality, that is, Young's inequality. If $u\in L^\Phi(\Omega)$ and $v\in L^\Psi(\Omega)$, with $(\Phi,\Psi)$ as a complementary Young pair, then we have $$\int_\Omega|uv|\,d\mu\leq 2\| u \|_{L^\Phi(\Omega)}\| v \|_{L^\Psi(\Omega)}.$$
Now we present a generalization of Young's inequality. Let $\Phi_i,\,i=1,2,3$, be Young functions for which $\Phi_1^{-1}(x)\Phi_2^{-1}(x)\leq \Phi_3^{-1}(x),\, x\geq 0$ holds. If $u_i\in L^{\Phi_i}(\Omega),\,i=1,2$, then $u_1\cdot u_2\in {L^{\Phi_3}(\Omega)}$ and 
$$\| u_1u_2\|_{L^{\Phi_3}(\Omega)}\leq 2\| u_1\|_{L^{\Phi_1}(\Omega)}\| u_2\|_{L^{\Phi_2}(\Omega)}.$$

\subsection{Metric-valued Orlicz-Sobolev spaces via upper gradients}\label{sec:metic valued via upper gradients}

Fix an N-function $\Phi$, we denote by $L^{\Phi}(\Omega,Y)$ the space of all $\mu$-measurable and essentially separably valued maps $u\colon \Omega\to Y$ such that for some $y_0\in Y$, the function $x\mapsto d(u(x),y_0)\in L^{\Phi}(\Omega)$. A sequence $\{u_k\}\subset L^{\Phi}(\Omega,Y)$ is said to converge to $u\in L^{\Phi}(\Omega,Y)$ if
$$\| d_Y(u,u_k)\|_{L^{\Phi}(\Omega)}\to 0 \qquad \text{as }k\to\infty.$$
When $(Y,d_Y)=(V, |\cdot|)$ is a Banach space, we may endow $L^{\Phi}(\Omega,V)$ with a natural norm 
$$\|f\|_{L^\Phi(\Omega, V)}:=\| |f| \|_{L^{\Phi}(\Omega)}.$$
If $V$ is $\real$, we simply set $L^{\Phi}(\Omega, \real):=L^{\Phi}(\Omega)$.

The concept of upper gradient was first introduced in~\cite{hk98,KM98} and was studied in details in \cite{s00}. We next introduce more general concept, namely, $\Phi$-weak upper gradients.  Here we only give a very brief introduction and refer the interested readers to the monograph \cite{H04} for more information.

\begin{definition}\label{def:upper gradients} 
	A Borel function $g\colon \Omega\rightarrow [0,\infty]$ is called an \emph{upper gradient} for a map $u\colon \Omega\to Y$ if for every rectifiable curve $\gamma\colon [a,b]\to \Omega$, we have the inequality
	\begin{equation}\label{upper-gradient}
		d_Y(u(\gamma(b)),u(\gamma(a)))\leq \int_\gamma g\,ds\text{.}
	\end{equation}
If inequality~\eqref{upper-gradient} holds for $\Phi$-almost every curve, then g is called a \emph{$\Phi$-weak upper gradient} for u.
\end{definition}
We say that a property of curves holds for $\Phi$-{\it almost every curve} if the collection of locally rectifiable curves for which the property fails to hold has $\Phi$-modulus zero. For definition and properties of $\Phi$-modulus, see \cite{H04}.

A $\Phi$-weak upper gradient $g$ of $u$ is \emph{minimal} if for every $\Phi$-weak upper gradient $\tilde{g}$ of $u$, $\tilde{g}\geq g$ $\mu$-almost everywhere.  If $u$ has an upper gradient in $L^{\Phi}_{\loc}(\Omega)$, then $u$ has a unique (up to sets of $\mu$-measure zero) minimal $\Phi$-weak upper gradient. We denote the minimal $\Phi$-weak upper gradient by $g_u$. The Orlicz-Sobolev space $N^{1,\Phi}(\Omega,Y)$ consists of all $u\in L^\Phi(\Omega,Y)$ with a minimal $\Phi$-weak upper gradient $g_u\in L^{\Phi}(\Omega)$. For each $u\in N^{1,\Phi}(\Omega,Y)$, we shall use $E^{\Phi}(u)$ to denote the $\Phi$-weak upper gradient energy functional of $u$, that is, 
$$E^{\Phi}(u)= \| g_u\|_{L^{\Phi}(\Omega)}.$$

An alternative way to introduce $N^{1,\Phi}(\Omega,Y)$ is to use isometric embedding $Y\subset V$ and then define $N^{1,\Phi}(\Omega,Y)$ as the Banach space-valued Sobolev spaces $N^{1,\Phi}(\Omega,V)$. We briefly record Banach space valued Sobolev spaces $N^{1,\Phi}(\Omega,V)$ here.

The {\it Dirichlet space} $D^{1, \Phi}(\Omega, V)$ consists of all measurable functions $u\colon \Omega\rightarrow V$ that have an upper gradient belonging to $L^\Phi(\Omega)$. We can equip the Dirichlet  space $D^{1, \Phi}(\Omega, V)$ with the seminorm
$$\|u\|_{D^{1, \Phi}(\Omega, V)}:=\inf_g\|g\|_{L^\Phi(\Omega)},$$
where the infimum is taken over all $\Phi$-weak upper gradient $g$ of $u$.

Let
$$\tilde N^{1,\Phi}(\Omega, V) =D^{1, \Phi}(\Omega, V) \cap L^\Phi(\Omega, V)$$
be equipped with the seminorm
$$\|u\|_{\tilde N^{1, \Phi}(\Omega, V)} =\|u\|_{L^\Phi(\Omega, V)}+\|u\|_{D^{1, \Phi}(\Omega, V)}.$$
We obtain a normed space $N^{1, \Phi}(X, V)$, which is called the {\it Sobolev space} of $V$-valued functions on $\Omega$, by passing to equivalence classes of functions in $\tilde N^{1, \Phi}(\Omega, V)$, where $u_1\sim u_2$ if and only if $\|u_1-u_2\|_{\tilde N^{1, \Phi}(\Omega, V)}=0$. Thus,  
$$N^{1, \Phi}(\Omega, V):=\tilde N^{1, \Phi}(\Omega, V)/\{u\in \tilde N^{1, \Phi}(\Omega, V): \|u\|_{\tilde N^{1, \Phi}(\Omega, V)=0}\}.$$
Since we may embed the metric space $Y$ isometrically into some Banach space   $L^\infty(Y)$, we can alternatively define $N^{1,\Phi}(\Omega,Z)$ via $N^{1,\Phi}(\Omega,Z):=N^{1,\Phi}(\Omega,L^\infty(Z))$.

A pair $(u,g)$, where u is locally integrable and $g\geq 0$, is said to satisfy the loacl $1$-{\it Poincar\'e inequality}, if there exist constants $C>0$ and $\lambda\geq 1$ such that
\begin{equation}\label{eq:poincare ineq}
	\vint_{B}|u-u_B|\,d\mu \leq Cr\vint_{\lambda B}g\,d\mu
\end{equation}
holds for each ball $B=B(x,r)$ satisfying $\lambda B\subset \Omega$. More generally, if $\Phi$  is a Young function, a pair as above is said to satiefy the local $\Phi$-{\it Poincar\'e inequality}, if there exist constants $C_\Phi>0$ and $\lambda\geq 1$ such that
\begin{equation}\label{eq:phi Poincare inequality}
	\vint_{B} |u-u_B|\, d\mu\leq C_\Phi r \Phi^{-1}\left(\vint_{\lambda B}\Phi(g)\,d\mu\right).
\end{equation}
holds for each ball $B=B(x,r)$ satisfying $\lambda B\subset \Omega$. The domain $\Omega$ supports the local $1$-Poincar\'e inequality (respectively $\Phi$-Poincar\'e inequality) if \eqref{eq:poincare ineq} (\eqref{eq:phi Poincare inequality}) holds for every integrable map $u$ and every upper gradient $g$ of $u$. Notice that \eqref{eq:poincare ineq} implies
\eqref{eq:phi Poincare inequality} by the Jensen inequality.


\subsection{Orlicz-Hajlasz-Sobolev spaces }
\begin{definition}[Hajlasz-Sobolev spaces]\label{def:Hajlasz-Sobolev}
	A measurable map $u\colon \Omega\to Y$ belongs to the \emph{Hajlasz-Sobolev space} $M^{1,\Phi}(\Omega,Y)$ if $u\in L^\Phi(X,Y)$ and there exists a nonnegative function $g\in L^\Phi(\Omega)$ such that the \emph{Hajlasz gradient inequality}
	\begin{equation}\label{eq:Hajlasz gradient}
		d_Y(u(x),u(x'))\leq d_X(x,x')\left(g(x)+g(x')\right)
	\end{equation} 
	holds for all $x,x'\in \Omega\backslash N$ for some $N\subset \Omega$ with $\mu(N)=0$. For each $u\in M^{1,\Phi}(\Omega,Y)$, the associated Hajlasz energy $E^\Phi_H(u)$ is defined as 
	$$E^\Phi_H(u):=\inf_{g}\|g\|_{L^{\Phi}(\Omega)},$$
	where the infimum is taken over all Hajlasz gradient $g$ of $u$, that is, $g$ such that \eqref{eq:Hajlasz gradient} holds.
\end{definition}

We need two maximal functions. Let $0\leq \alpha<\infty$, $0<\beta<\infty$ and $u\in L_{\loc}^1(\Omega)$. The \emph{fractional maximal function} of $u$ is defined to be
\begin{equation}\label{eq:maximal operator}
	M_{\alpha}u(x)=\sup\limits_{r>0}r^\alpha\vint_{B(x,r)}|u|\,d\mu,
\end{equation}
where $u_B=\vint_B u\,d\mu={\mu(B)}^{-1}\int_B u\,d\mu$ is the integral average of $u$ over $B$. If $\alpha =0$, then we obtain the classical \emph{Hardy-Littlewood maximal function}. As is well known,  $M_0\colon L^\Phi(X)\to L^\Phi(X)$ is bounded if and only if $\Phi\in\nabla_2$, see \cite{C99}. The \emph{fractional sharp maximal function} of $u$ is defined as
\begin{equation}\label{eq:sharp maximal operator}
	u_\beta^\sharp(x)=\sup\limits_{r>0}r^{-\beta}\vint_{B(x,r)}|u-u_{B(x,r)}|\,d\mu.
\end{equation}

Next we give two sufficient conditions for $\Omega$ to be admissible.

\begin{lemma}\label{lem:m=n by 1 poincare}
	Assume that $\Phi$ is a doubling N-function, the maximal operator $M_0$ defined in \eqref{eq:maximal operator} is bounded in $L^\Phi(\Omega)$ and that  $\Omega$ supports a  $1$-Poincar\'e inequality, then $N^{1,\Phi}(\Omega, Y) = M^{1,\Phi}(\Omega, Y)$ with comparable norms for all Orlicz-Sobolev maps.
\end{lemma}

\begin{proof}
	In general, $M^{1,\Phi}(\Omega, Y) \subset N^{1,\Phi}(\Omega, Y)$. For the other direction, let $u\in N^{1,\Phi}(\Omega, Y)$ with a $\Phi$-weak upper gradient $g_u\in L^{\Phi}(\Omega)$. By \cite[Proof of Lemma 3.6]{Hk98}, the inequality $$d_Y(u(x),u(y))\leq C d(x,y)(u_1^\sharp(x)+u_1^\sharp(y))$$
	holds for almost every $x,y\in\Omega$. Note that the $1$-Poincar\'e inequality implies that $$u_1^\sharp(x)\leq CM_0 g_u(x)$$ for all $x\in\Omega$. The claim follows from the boundness of $M_0$.

\end{proof}

We shall need the following result for the next lemma.

\begin{proposition}[\cite{H04}, Lemma 5.15]\label{pro:maximal operator ineq}
	Assume that $\Phi\colon\left[0,\infty\right)\to\left[0,\infty\right)$ is a strictly increasing Young function. If a pair $u\in L_{\loc}^1(\Omega,Y)$ and a measurable function $g\geq 0$ satisfy a $\Phi$-Poincar\'e inequality, then for $\mu$-almost all $x,y\in \Omega$,
	$$d_Y(u(x),u(y))\leq Cd(x,y)\left(\Phi^{-1}\left(M_0\Phi(g(x))\right)+\Phi^{-1}\left(M_0\Phi(g(y))\right)\right),$$
	where the constant $C>0$ depends only on the doubling constant $C_d$ of $\mu$ and on the constant $C_\Phi$ of the $\Phi$-Poincar\'e inequality.
\end{proposition}

\begin{lemma}\label{lem:m=n by phi poincare}
	Assume that $\Phi$ and $\Phi_1$ are N-functions and $\Phi_1$ belongs to $\nabla_2$. If the function $\Phi_2=\Phi_1^{-1}\circ\Phi$ is a doubling N-function and if $\Omega$ supports a $\Phi_2$-Poincar\'e inequality, then $N^{1,\Phi}(\Omega, Y) = M^{1,\Phi}(\Omega, Y)$ with comparable norms for all Orlicz-Sobolev maps.
\end{lemma}

\begin{proof}
	In general, $M^{1,\Phi}(\Omega, Y) \subset N^{1,\Phi}(\Omega, Y)$. Therefore, we just need to show that $N^{1,\Phi}(\Omega, Y) \subset M^{1,\Phi}(\Omega, Y)$. To show this, let $u\in N^{1,\Phi}(\Omega, Y)$ with an $\Phi$-weak upper gradient $g_u\in L^{\Phi}(\Omega)$. By Proposition~\ref{pro:maximal operator ineq}, 
	$$d_Y(u(x),u(y))\leq Cd(x,y)\left(\Phi_2^{-1}\left(M_0\Phi_2(g_u(x))\right)+\Phi_2^{-1}\left(M_0\Phi_2(g_u(y))\right)\right)$$
	for $\mu$-almost all $x,y\in \Omega$. It suffices to show that the function $\Phi_2^{-1}(M_0\Phi_2(g_u))$ belongs to $L^\Phi(\Omega)$. Since $g_u$ belongs to $L^\Phi(\Omega)$, $\Phi_2(g_u)$ is in $L^{\Phi_1}(\Omega)$. As the maximal operator is bounded in $L^{\Phi_1}$ by the $\nabla_2$-property of $\Phi_1$, the function $M_0\Phi_2(g_u)$ belongs to $L^{\Phi_1}(\Omega)$. Consequently, by the definition of $\Phi_2$, $\Phi_2^{-1}(M_0\Phi_2(g_u))$ is a function of $L^\Phi(\Omega)$. Hence $u$ belongs to $M^{1,\Phi}(\Omega, Y)$. The proof is complete.
\end{proof}

\subsection{Orlicz-Sobolev embedding theorem}
In the classical Sobolev-Poincar\'e inequalities in $\rr^n$, the dimension $n$ plays a special role, particularly in the embedding theorems of Sobolev and Morrey. For our visions of the Sobolev-Poincar\'e-type estimate, a suitable substitute for this threshold parameter is given by a lower decay order for the measure of balls. We say that X has \emph{a relative lower volume decay of order $s>0$} if
\begin{equation}\label{eq:relative lower volume decay of order}
	\left(\frac{\diam(B')}{\diam(B)}\right)^s\leq C_0\frac{\mu(B')}{\mu(B)}
\end{equation} 
holds whenever $B'\subset B$ are balls in X. Note $\mu$ is doubling, inequality \eqref{eq:relative lower volume decay of order} always holds for some $s\leq\log_2 C_d$ whenever $B'\subset B$ are balls in X, where $C_d$ is doubling constant of $\mu$.


Let $s>1$. For a Young function $\Phi$ satisfying
\begin{equation}\label{eq:tiaojian1}
	\int_0\left(\frac{t}{\Phi(t)}\right)^{s'-1}\,dt<\infty\quad \text{and}\quad \int^\infty\left(\frac{t}{\Phi(t)}\right)^{s'-1}\,dt=\infty,
\end{equation}
we define
\begin{equation}\label{eq:tiaojian2}
	\Phi_s=\Phi\circ\Psi_s^{-1},
\end{equation}
where
\begin{equation}\label{eq:tiaojian3}
	\Psi_s(r)=\left(\int_0^r\left(\frac{t}{\Phi(t)}\right)^{s'-1}\,dt\right)^{1/s'} \text{and}\quad s'=\frac{s}{s-1}.
\end{equation}
Here, $\int_0 f(t)\,dt<\infty$ means that $\int_0^c f(t)\,dt<\infty$ for some $c>0$. Similarly, we denote $\int^\infty f(t)\,dt=\infty$, if $\int_c^\infty f(t)\,dt=\infty$ for all $c>0$. In the classical Euclidean case, these formulas
\eqref{eq:tiaojian1}-\eqref{eq:tiaojian3} were first introduced in \cite{C97}.

The following Orlicz-Sobolev embedding theorem is well-known.
\begin{proposition}[\cite{HK22}, \cite{H10}]\label{thm:embedding thm}
	 Assume that $(X,d,\mu)$ is a doubling metric measure space that supports the $\Phi$-Poincar\'e inequality \eqref{eq:phi Poincare inequality} and satisfies \eqref{eq:relative lower volume decay of order} with $s>1$. Let $B$ be a ball, $\delta>0$ and $\tilde{B}=(1+\delta)\lambda B$. If $\Phi$ satisfies \eqref{eq:tiaojian1}, 
	then $N^{1,\Phi}(\tilde B)\subset L^{\Phi_s}(B)$, where $\Phi_s$ is defined by \eqref{eq:tiaojian2} and \eqref{eq:tiaojian3}. 
	
	Moreover, for every $u\in N^{1,\Phi}(\tilde B)$ and for every $\Phi$-weak upper gradient $g$ of $u$, we have 
	 $$\| u-u_B\|_{L^{\Phi_s}(B)}
	 \leq C r_B \mu(B)^{-1/s}\| g\|_{L^{\Phi}(\tilde B)}.$$
\end{proposition}

\section{Generalized Rellich-Kondrachov theorem}\label{sec:proof}

\subsection{Proof of Theorem \ref{thm:main theorem}}

We firstly prove a Rellich-Kondrachov compactness theorem for Orlicz-Hajlasz-Sobolev mappings.
\begin{theorem}\label{thm:proof theorem}
	Let $\mu$ be a doubling measure on $X$ and $\Omega$ supports a loacl $\Phi$-Poincar\'e inequality for some  N-function $\Phi\in\Delta'\cap\nabla_2$. For every $m\in\nn$, let $(Y_m,d_m)$ be a complete metric space, $K_m\subset Y_m$ compact and $\{u_m\} \subset M^{1,\Phi}(\Omega, Y_m)$. Suppose that $(K_m, d_m)$ is a uniformly compact sequence and 
	\begin{equation}\label{eq:bdd-energy-distance-seq-um}
		\sup_{m\in\nn}\left[ \| d_m(y_m,u_m)\|_{L^\Phi(\Omega)}+E_H^\Phi(u_m)\right]<\infty
	\end{equation}
	for some and thus every $y_m\in K_m$. Then, after possibly passing to a subsequence, there exist a complete metric space $Z$, isometric embeddings $\varphi_m: Y_m \hookrightarrow Z$, a compact subset $K\subset Z$ and $v\in M^{1,\Phi}(\Omega, Z)$ such that $\varphi_m(K_m)\subset K$ for all $m\in \nn$ and $\varphi_m\circ u_m$ converges to $v$ in $L^\Phi(\Omega, Z)$.
\end{theorem}

\begin{proof}[Proof of Theorem~\ref{thm:proof theorem}]

As in \cite{GW20}, the proof of theorem~\ref{thm:proof theorem} relies on the following variant of Gromov's compactness theorem.

\begin{proposition}\label{prop:extended-gromov-compactness-sets}
Let  $(Y_m, d_m)$ be a sequence of metric spaces and, for each $m\in\nn$, subsets
\begin{equation*}
	B_m^1\subset B_m^2\subset B_m^3\subset\dots\subset Y_m.
\end{equation*}
If for every $k\in\nn$ the sequence $(B_m^k, d_m)$ is uniformly compact then, after possibly passing to a subsequence, there exist a complete metric space $Z$, isometric embeddings $\varphi_m\colon Y_m\hookrightarrow Z$ and compact subsets $Y^1\subset Y^2\subset\dots\subset Z$ such that $\varphi_m(B_m^k)\subset Y^k$ for all $m\in\nn$ and $k\in\nn$.
\end{proposition}

Fix $m\in\nn$. Reasoning as in \cite[Proof of Theorme 3.1]{GW20}, we may assume $Y_m$ is a Banach space. Indeed, every metric space $Y$ isometrically embeds into the Banach space $L^\infty(Y)$ of bounded functions on $Y$ with the supremum norm. Now, there exists a non-negative function $h_m\in L^\Phi(\Omega)$ such that $\|h_m\|_{L^\Phi(\Omega)}\leq C\cdot E_H^\Phi(u_m)$ for some constant $C$ only depending on $\Omega$ and $\Phi$ and such that
\begin{equation*}
	d_m(u_m(x), u_m(x')) \leq d(x, x') (h_m(x) + h_m(x'))
\end{equation*}
for all $x,x'\in \Omega$. For $k\in \nn$ set $$A_m^k:=\{ x\in \Omega : h_m(x) \leq k\}$$ and notice that the restriction of $u_m$ to $A_m^{k}$ is $2k$-Lipschitz.

\begin{lemma}\label{lem:bdd-diam}
There exist $k_0\in\nn$ and $\lambda>0$ such that $u_m(A_m^k)\subset B(K_m, \lambda k)$ and $A_m^{k}\not=\emptyset$ for all $m\in\nn$ and $k\geq k_0$.
\end{lemma}

Here, $B(K_m, \lambda k)$ denotes the set of all $y\in Y_m$ for which there exists $z\in K_m$ with $d_m(y,z)<\lambda k$.

\begin{proof}
	For each $m\in\nn$, fix $y_m\in K_m$ and define 
	$$C_m^k:= \{x\in\Omega : d_m(y_m, u_m(x))\leq k\}.$$
	By $\Phi\in\Delta'$ and
	\eqref{eq:bdd-energy-distance-seq-um}, there exists $M>0$ such that
	\begin{equation*}\label{eq:bdd-Amk}
		\begin{split}
		\mu(\Omega\setminus A_m^k)
		&\leq \Phi(1)^{-1} \int_{\Omega\setminus A_m^k }\Phi\left(\frac{h_m(x)}{k}\right)\,d\mu\\
		&\leq M'\Phi(1)^{-1} \Phi\left(\frac{1}{k}\right) \int_\Omega \Phi(h_m(x))\,d\mu\\
		&\leq M\Phi\left(\frac{1}{k}\right)
		\end{split}
\end{equation*} 
and similarly we have $\mu(\Omega\setminus C_m^k)\leq M\Phi\left(\frac{1}{k}\right)$ for all $m$ and $k$. Thus, there exists $k_0\in\nn$ such that $A_m^k\cap C_m^k\not=\emptyset$ for all $m\in\nn$ and all $k\geq k_0$. Fix $x_0\in A_m^k\cap C_m^k$. Then for every $x\in A_m^k$ we have 
$$d_m(y_m, u_m(x)) \leq d_m(y_m, u_m(x_0)) + d_m(u_m(x_0), u_m(x)) \leq k+ 2k\diam(\Omega),$$ so the lemma follows. 
\end{proof}

Let $m\in\nn$ and $k\geq k_0$. By \cite[Theorem~4.1.21]{hkst12}, there exists a $2kC_d$-Lipschitz map $u_m^k\colon\overline{\Omega}\to Y_m$ which agrees with $u_m$ on $A_m^k$, where $C_d$ is the doubling constant of $\mu$. We define for each $m\in\nn$ an increasing sequence of subsets $B_m^{k_0}\subset B_m^{k_0+1}\subset \dots \subset Y_m$ by $$B_m^k:= K_m\cup u_m^{k_0}(\overline{\Omega})\cup\dots\cup u_m^k(\overline{\Omega}).$$ Since $u_m^j$ is $2jC_d$-Lipschitz on the compact set $\overline{\Omega}$, Lemma~\ref{lem:bdd-diam} implies that for fixed $k\geq k_0$ the sequence of metric spaces $(B_m^k, d_m)$ is uniformly compact.
Thus, by Proposition~\ref{prop:extended-gromov-compactness-sets} there exists, after possibly passing to a subsequence, a complete metric space $(Z, d_Z)$, isometric embeddings $\varphi_m\colon Y_m\hookrightarrow Z$, and compact subsets $Y^{k_0}\subset Y^{k_0+1}\subset\dots\subset Z$ such that $\varphi_m(B_m^k)\subset Y^k$ for all $m$ and $k\geq k_0$. In particular, for every $m\in\nn$ the set $\varphi_m(K_m)$ is contained in the compact set $K:=Y^{k_0}$. 

\begin{lemma}\label{lemma:vm-w}
	The maps $v_m = \varphi _m\circ u_m$ belong  to $M^{1,\Phi}(\Omega, Z)$ and satisfy
	\begin{equation}\label{eq:bounded-vj-inZ}
		\sup_{m\in\nn}\left[\| d_Z(z_0,v_m)\|_{L^{\Phi}(\Omega)}+ E_H^\Phi(v_m)\right]<\infty
	\end{equation}
	for some and thus every $z_0\in Z$.
\end{lemma}

\begin{proof}
	Fix $x'\in\Omega$ and write $z_0=v_m(x')$. Then we have
	\begin{equation*}
		\begin{split}
		d_Z(z_0,v_m(x))&= d_Z(v_m(x'),v_m(x))\\
		&= d_Z(\varphi_m\circ u_m(x'),\varphi_m\circ u_m(x))
		=d_m(u_m(x'),u_m(x)).
		\end{split}
	\end{equation*}
Since $u_m \in M^{1,\Phi}(\Omega, Y_m)$, then $d_m(u_m(x'),u_m)\in L^{\Phi}(\Omega)$. Hence $d_Z(z_0,v_m)\in L^{\Phi}(\Omega)$.
For the Hajlasz gradient inequality, note that
\begin{equation*}
	\begin{split}
		d_Z(v_m(x'),v_m(x))&= d_Z(\varphi_m\circ u_m(x'),\varphi_m\circ u_m(x))\\
		&=d_m(u_m(x'),u_m(x))
		\leq d(x',x)(h_m(x')+h_m(x))
	\end{split}
\end{equation*}
for all $x, x'\in\Omega$. Therefore, $v_m\in
M^{1,\Phi}(\Omega, Z)$ and \eqref{eq:bounded-vj-inZ} follows form \eqref{eq:bdd-energy-distance-seq-um}.
\end{proof}

We need the following lemma for the later proofs. 

\begin{lemma}\label{lem:diff-um-umk}
	For given $k\geq k_0$, define the map $v_m^k:= \varphi_m\circ u_m^k$. Then there exists $M>0$ such that 
	\begin{equation}
		\int_\Omega d_Z(v_m(x), v_m^k(x))\,d\mu \leq M\left[k\Phi\left(\frac{1}{k}\right)+\Psi^{-1}\left(M\Phi\left(\frac{1}{k}\right)\right)\right]
	\end{equation}
	for all $m\in\nn$ and every $k\geq k_0$, where $\Psi$ is the complementary function of $\Phi$. 
\end{lemma}

\begin{proof}
     Fix $x'\in \Omega$ and write $z_0=v_m^k(x')$. Then 
	\begin{equation*}
		\begin{split}
			d_Z(z_0,v_m^k(x))&= d_Z(\varphi_m\circ u_m^k(x'),\varphi_m\circ u_m^k(x))\\
			&= d_m(u_m^k(x'), u_m^k(x))
			\leq 2kC_d \diam(\Omega)
		\end{split}
	\end{equation*}
	for every $x\in\Omega$ and all $m\in\nn$ and $k\geq k_0$. This together with Young's inequality and Lemma~\ref{lem:bdd-diam} yields
	\begin{equation*}
		\begin{split}
			\int_\Omega d_Z(v_m(x), v_m^k(x))\,d\mu &= \int_{\Omega\setminus A_m^k} d_Z(v_m(x), v_m^k(x))\,d\mu\\
			&= \int_{\Omega\setminus A_m^k} d_Z(v_m^k(x), z_0)\,d\mu+\int_{\Omega\setminus A_m^k} d_Z(v_m(x), z_0)\,d\mu\\
			&\leq M'k\mu(\Omega\setminus A_m^k)+2\| d_Z(v_m(x),z_0)\|_{L^{\Phi}(\Omega\setminus A_m^k)}\| 1 \|_{L^{\Psi}(\Omega\setminus A_m^k)}\\
			&\leq M'k\mu(\Omega\setminus A_m^k)+M'''\Psi^{-1}(M''\mu(\Omega\setminus A_m^k))\\
			&\leq M\left[k\Phi\left(\frac{1}{k}\right)+\Psi^{-1}\left(M\Phi\left(\frac{1}{k}\right)\right)\right],
		\end{split}
	\end{equation*}
	where $(\Phi,\Psi)$ is a complementary Young pair, the constants $M,\,M',\,M''$ and $M'''$ do not depend on $m$ and $k$.
\end{proof}

We next derive the $L^1$ limit.
\begin{lemma} \label{lem:vmj-conv-L1}
There exists a subsequence $\{v_{m_j}\}$ which converges in $L^1(\Omega, Z)$ to some $v\in L^1(\Omega, Z)$.
\end{lemma}

\begin{proof}
	For given $k\geq k_0$, the map $v_m^k:= \varphi_m\circ u_m^k$ is $2kC_d$-Lipschitz and has image in the compact set $Y^k$ for every $m\in\nn$. Thus, by the Arzel\`a-Ascoli theorem and by a diagonal sequence argument, there exist  integers $1\leq m_1<m_2<\dots$ such that, for every $k\geq k_0$, the sequence $\{v_{m_j}^k\}$ converges uniformly on $\Omega$ as $j\to\infty$. Lemma~\ref{lem:diff-um-umk} above shows that there exists $M>0$ such that
	\begin{equation*}
		\int_\Omega d_Z(v_{m_j}(x), v_{m_l}(x))\,d\mu \leq 2M\left[k\Phi\left(\frac{1}{k}\right)+\Psi^{-1}\left(M\Phi\left(\frac{1}{k}\right)\right)\right] + \int_\Omega d_Z(v_{m_j}^k(x), v_{m_l}^k(x))\,d\mu
	\end{equation*}
	for all $j, l\in\nn$ and every $k\geq k_0$. Hence, the integral on the left-hand side converges to $0$ as $j,l\to\infty$. This proves that $\{v_{m_j}\}$ is a Cauchy sequence in $L^1(\Omega, Z)$ and hence that $v_{m_j}$ converges in $L^1(\Omega, Z)$ to some $v\in L^1(\Omega, Z)$.
\end{proof}

\begin{lemma}\label{lem:conv-in-Lp}
	The sequence $\{v_{m_j}\}$ converges to $v$ in $L^{\Phi}(\Omega, Z)$.
\end{lemma}
\begin{proof}
	By Lemma~\ref{lemma:vm-w}, there exists some $v'\in L^{\Phi}(\Omega,Z)$ such that $v_{m_j}\rightharpoonup v'$ in $L^{\Phi}(\Omega,Z)$. Since $v_{m_j}\to v$ in $L^1(\Omega,Z)$, $v=v'$ and so $v\in L^{\Phi}(\Omega,Z)$.
	
	Let $\varepsilon>0$. Then the set $F_\varepsilon^j:= \{x\in \Omega: d_Z(v_{m_j}(x), v(x)) >\varepsilon\}$ satisfies $\mu (F_\varepsilon^j)\to 0$ as $j\to\infty$ because, by Chebyshev's inequality, $$\mu(F_\varepsilon^j) \leq \varepsilon^{-1} \cdot \int_\Omega d_Z(v_{m_j}(x), v(x))\,d\mu$$ for every $j\in\nn$ and because $v_{m_j}$ converges to $v$ in $L^1(\Omega, Z)$ by Lemma~\ref{lem:vmj-conv-L1}.
	
	By the absolute continuity of the integral, we have
	\begin{equation*}
		\begin{split}
			\int_\Omega\Phi(d_Z(v_{m_j}(x),v(x)))\,d\mu&\leq \int_{\Omega\setminus F_\varepsilon^j}\Phi(d_Z(v_{m_j}(x),v(x)))\,d\mu+\int_{F_\varepsilon^j} \Phi(d_Z(v_{m_j}(x),v(x)))\,d\mu\\
			&\leq \int_{\Omega\setminus F_\varepsilon^j}\Phi(\varepsilon)\,d\mu+\int_{F_\varepsilon^j} \Phi(d_Z(v_{m_j}(x),v(x)))\,d\mu
		\end{split}
	\end{equation*}
		and hence 
		$$\| d_Z(v_{m_j},v)\|_{L^{\Phi}(\Omega)}\to 0$$
		as $j\to \infty$. Above, we used the well-known fact that if $\Phi$ is doubling, then $L^{\Phi}$-convergence is equivalent to $\Phi$-mean convergence; see \cite[Theorem 3.4.12]{RR91}. This shows that $v_{m_j}$ converges to $v$ in $L^\Phi(\Omega, Z)$.
\end{proof}

Finally, we are able to finish the proof of Theorem \ref{thm:proof theorem}.
\begin{lemma}
	Suppose that $\{v_{m_j}\}\in M^{1, \Phi}(\Omega, Z)$ and $v_{m_j}\to v$ in $L^1(\Omega,Z)$. Then $v\in M^{1,\Phi}(\Omega, Z)$.
\end{lemma}
\begin{proof}
	Since $\{v_{m_j}\}\in M^{1, \Phi}(\Omega, Z)\}$, then there exist nonnegative functions $\{h_{m_j}\}\subset L^\Phi(\Omega
	)$ such that $$d_Z(v_{m_j}(x),v_{m_j}(y))\leq d_X(x,y)\left(h_{m_j}(x)+h_{m_j}(y)\right)$$
	for all $x,y\in \Omega$.
	Since $E_H^\Phi(v_{m_j})$ is bounded, then, after possibly passing to a subsequence, there exists a function $h\in L^\Phi(\Omega)$ such that $h_{m_j}\rightharpoonup h$ in $L^\Phi(\Omega)$. Since $v_{m_j}\to v$ in $L^1(\Omega,Z)$, then, up to a subsequence, we have $v_{m_j}\to v$ alomst everywhere in $\Omega$. Similarly, $h_{m_j}\to h$ alomst everywhere in $\Omega$. Hence, we have
	$$d_Z(v(x),v(y))\leq d_X(x,y)\left(h(x)+h(y)\right)$$
	for almost all $x,y\in \Omega$. This together with Lemma~\ref{lem:conv-in-Lp} yields $v\in M^{1,\Phi}(\Omega, Z)$, completing the proof.
\end{proof}

The proof of Theorem \ref{thm:proof theorem} is thus complete.
\end{proof}

\begin{proof}[Proof of Theorem~\ref{thm:main theorem}]
     With Lemma~\ref{lem:m=n by 1 poincare} and Lemma~\ref{lem:m=n by phi poincare} at hand, Theorem \ref{thm:main theorem} follows immediately from Theorem \ref{thm:proof theorem}. 
\end{proof}

\subsection{Some further improvements}

In Theorem \ref{thm:main theorem}, we proved that $\varphi_m\circ u_m\to v$ in $L^{\Phi}$ and this is to some extent weaker than the statement of classical Rellich-Kondrachov compactness theorem, as we did not get the convergence in a Orlicz space better than $L^{\Phi}$. In this subsection, we discuss how to improve the $L^{\Phi}$ convergence. 


\begin{corollary}\label{lem:conv-in-Lq}
	For every $N$-function $H\prec\prec\Phi_s$ the maps $\{v_{m_j}\}$ and $v$ belong to $L^H(\Omega, Z)$ and the sequence $(v_{m_j})$ converges to $v$ in $L^H(\Omega, Z)$, where $\Phi_s$ is defined by \eqref{eq:tiaojian2} and \eqref{eq:tiaojian3}.
\end{corollary}

\begin{proof} 
	Fix $z_0\in Z$. If $H\prec\prec\Phi$, then Corollary~\ref{lem:conv-in-Lq} follows from Theorem \ref{thm:main theorem}. Therefore, we may assume $\Phi\prec\prec H\prec\prec \Phi_s$. Note that $M^{1,\Phi}(\Omega, Z)\subset N^{1,\Phi}(\Omega, Z)$, by \eqref{eq:bounded-vj-inZ} and the Orlicz-Sobolev embedding (see Proposition~\ref{thm:embedding thm}), the real-valued functions $x\mapsto d_Z(z_0, v_{m_j}(x))$ belong to $L^{\Phi_s}(\Omega)$ and form a bounded sequence in $L^{\Phi_s}(\Omega)$. 
	
	Since a subsequence of $\{v_{m_j}\}$ converges to $v$ almost everywhere, it follows with Fatou's lemma that $v\in L^{\Phi_s}(\Omega, Z)$ and hence $$L:= \sup_{j\in\nn} \| d_Z(v_{m_j},v)\|_{L^{\Phi_s}(\Omega)}
	<\infty.$$  
	By the linearity of the integral, we have
	\begin{equation*}
		\begin{split}
			\int_\Omega H(d_Z(v_{m_j}(x),v(x)))\,d\mu&\leq \int_{\Omega\setminus F_\varepsilon^j}H(d_Z(v_{m_j}(x),v(x)))\,d\mu+\int_{F_\varepsilon^j}H(d_Z(v_{m_j}(x),v(x)))\,d\mu\\
			&\leq \int_{\Omega\setminus F_\varepsilon^j}H(\varepsilon)\,d\mu+\int_{F_\varepsilon^j} H(d_Z(v_{m_j}(x),v(x)))\,d\mu
		\end{split}
	\end{equation*} 
	This together with generalized Young's inequality yields 
	\begin{align*}
		\| d_Z(v_{m_j}(x),v(x))\|_{L^{H}(\Omega)}&\lesssim \| \varepsilon \|_{L^{H}(\Omega)}+\| d_Z(v_{m_j}(x),v(x))\|_{L^{H}(\F_\varepsilon^j)}\\
		&\lesssim \| \varepsilon \|_{L^{H}(\Omega)}+\| d_Z(v_{m_j}(x),v(x))\|_{L^{\Phi_s}(\F_\varepsilon^j)}\| 1 \|_{L^{ P}(\F_\varepsilon^j)} 
	\end{align*}
	where $P$ is any Young function satisfying $\Phi_s^{-1}(t)P^{-1}(t)\leq H^{-1}(t)$ for all $t\geq 0$. Indeed, here we can simply set $P^{-1}(t):=H^{-1}\left(\frac{Ct}{H\circ \Phi^{-1}(t)}\right)$ and it follows easily from the property $H\prec\prec \Phi_s$ such that $P$ is a desired Young function with the above property. Since $\mu(F_\varepsilon^j)\to 0$ as $j\to\infty$, $\| 1 \|_{L^{ P}(\F_\varepsilon^j)}\to 0$ as $j\to \infty$. Consequently, 
	\begin{equation*}
		\| d_Z(v_{m_j}(x),v(x))\|_{L^{H}(\Omega)}\to 0
	\end{equation*}
	as $j\to \infty$. This shows that $v_{m_j}$ converges to $v$ in $L^H(\Omega, Z)$. The proof is complete.
\end{proof}

\section{Extension of the trace theory of Korevaar-Schoen}\label{sec:trace}

Let $(X, d_X,\mu)$ be a complete metric measure space, $(V, |\cdot|)$ be a Banach space  and $\Omega\subset X$ be a bounded domain. Assume $\bx$ is endowed with  an upper codimension-$\theta$ regular measure $\mathcal H$ with $\theta>0$. 

We give an alternative definition of the trace for Banach valued maps.
\begin{definition}\label{trace-Banach}
	Let $u\colon \Omega\rightarrow V$ be a $\mu$-measurable map. Then $Tu(x)\in V$ is the trace of $u$ at $x\in \bx$ if the following equation holds:
	\begin{equation}\label{trace-defn}
		\lim_{r\rightarrow 0^+}\vint_{B(x, r)\cap\Omega}|u-Tu(x)|\, d\mu=0.
	\end{equation}
	We say that $u$ has a trace $Tu$ on  $\bx$ if $Tu(x)$ exists for $\mH$-almost every $x\in \bx$.
	
	Arguing as in \cite[Lemma 3.2]{GHWX21}, we know that Definition \ref{trace-Banach} is consistent with Definition \ref{trace-metric}. Thus to develop a theory of trace, we shall not distinguish the trace operators in Definitions \ref{trace-metric} and \ref{trace-Banach}. From now on, we shall focus on the case when $(Y,d_Y)=(V,|\cdot|)$ is a Banach space.
	
\end{definition}
For any $f\in L^\Phi_{\loc}(\Omega, V)$, we define the centered fractional maximal operator as
\begin{equation}\label{maximal-operator}
	M_{\theta, \Phi} f(x)=\sup_{0<r<2\diam (\bx)} \Phi^{-1}\left(r^\theta\vint_{B(x, r)\cap \Omega} \Phi(f)\, d\mu\right),\ \ \text{ for each } \ x\in \bx.
\end{equation}
Then it is easy to see that this fractional maximal operator maps $L^\Phi_{\rm loc}(\Omega, V)$ into the space of real-valued lower semicontinuous functions on $\bx$.

Let us now establish its boundedness. The essential idea of the proof is similar with the one used in \cite[Lemma 4.2]{m17}.

\begin{lemma}\label{fractional max operator}
	Let the N-function $\Phi$ belong to $\Delta'$. Then for any $f\in L^\Phi(\Omega)$ and $\|f\|_{L^\Phi(\Omega)}\leq1$, we have $\|M_{\theta,\Phi}f\|_{WL^\Phi(\bx)}\leq\|f\|_{L^\Phi(\Omega)}$. In particular, $\mH(\{x\in\bx: M_{\theta,\Phi}f(x)=\infty\})=0$.
\end{lemma}

Here, given a Young function $\Phi$, the weak Orlicz space
$$WL^\Phi(\Omega)=\left\{f\,\text{is}\, \mu\text{-measurable}\colon\|f\|_{WL^\Phi(\Omega)}<\infty\right\},$$
where
$\|f\|_{WL^\Phi(\Omega)}=\sup\limits_{t>0}\displaystyle\int_\Omega \Phi(t\chi_{(t,\infty)}(|f|))\,d\mu.$

\begin{proof}
	Let $f\in L^\Phi(\Omega)$. Fix $t>0$, we can define $E_t=\{x\in\bx:M_{\theta,\Phi}f(x)>t\}$. For each $x\in E_t$, there is a ball $B_x=B(x,r_x)$ such that $$r_x^\theta\vint_{B_x\cap \Omega}\Phi(f)\,d\mu>\Phi(t
	).$$ Hence, $E_t\subset\bigcup_{x\in E_t}B_x\cap\bx$. Since $\mu|_\Omega$ is doubling and radii $r_x$ are bounded by $2\diam(\bx)$, we can apply the Basic Covering Theorem to find pairwise disjoint balls $B_k:=B_{x_k}$, $k=1,2,\cdots$, for some choice of $\{x_k\}\subset E_t$ such that $E_t\subset\bigcup_k 5B_k\cap\bx$. This together with upper codimension relation \eqref{upper} yields 
	$$\mH(E_t)\leq\mH \left(\bigcup\limits_{k}5B_k\cap\bx\right)\lesssim\sum\limits_{k}\mH(5B_k\cap\bx)\lesssim\sum\limits_{k}\frac{\mu(5B_k\cap\Omega)}{(5r_k)^\theta}\approx \sum\limits_{k}\frac{\mu(B_k\cap\Omega)}{r_k^\theta}.$$
	By the choice of balls $B_x$, we have
	$$\frac{\int_{B_k}\Phi(f)\,d\mu}{\Phi(t)}>\frac{\mu(B_k\cap\Omega)}{r_k^\theta}.$$
    That gives us that
    $$\mH(E_t)\lesssim\sum\limits_{k}\left(\frac{\int_{B_k}\Phi(f)\,d\mu}{\Phi(t)}\right)\lesssim \frac{\sum_k\int_{B_k}\Phi(f)\,d\mu}{\Phi(t)}\leq \frac{\int_\Omega\Phi(f)\,d\mu}{\Phi(t)}.$$
    It follows from the $N$-function $\Phi\in\Delta'$ 
    that
    \begin{align*}
    	\|M_{\theta,\Phi}f\|_{WL^{\Phi}(\bx)} 
    	&=\sup\limits_{t>0}\int_{\bx}\Phi(t\chi_{(t,\infty)}(|M_{\theta,\Phi}f(x)|))\,d\mH(x)\\
    	&\lesssim\sup\limits_{t>0}\left(\Phi(t)\int_{\bx}\Phi(\chi_{(t,\infty)}(|M_{\theta,\Phi}f(x)|))\,d\mH(x)\right)\\
    	&\leq\sup\limits_{t>0}\left(\Phi(t)\int_{E_t}\Phi(1)\,d\mH(x)\right)\\
    	&\lesssim\sup\limits_{t>0}\left(\Phi(t)\mH(E_t)\right)\\
    	&\leq\int_\Omega\Phi(f)\,d\mu\leq\|f\|_{L^{\Phi}(\Omega)}.
    \end{align*}
    	For the second claim, we assume that $\mH(\{x\in\bx: M_{\theta,\Phi}f(x)=\infty\})=a>0$. Set $$E_\infty=\{x\in\bx: M_{\theta,\Phi}f(x)=\infty\}.$$
    	For any $M>0$, we have $\|M_{\theta,\Phi}f\|_{WL^{\Phi}(\bx)}\geq\int_{E_\infty}\Phi(M)\,d\mH$, this is a contradiction as $M\to \infty$. So the lemma follows.
    
\end{proof}

We are ready to prove the boundedness of the trace operator for Banach space valued Orlicz-Sobolev maps. When $\Phi(x)=|x|^p$, $1<p<\infty$, the result was obtained in \cite[Theorem 3.4]{GHWX21}. The essential idea of the proof is similar with the one used in \cite{GHWX21}.

\begin{theorem}\label{thm:trace bdd}
	Suppose $\Omega\subset X$  is weakly $(\Phi,\theta)$-admissible for some N-function $\Phi\in\Delta'\cap\nabla_2$. 
	Then the trace operator $T\colon N^{1, \Phi}(\Omega, V)\rightarrow L^\Phi(\bx, V)$ is bounded and linear.
\end{theorem}

\begin{proof}
	Let $u\in N^{1, \Phi}(\Omega, V)$ and $R=2\diam(\Omega)$ be fixed. We extend $u$ as $0$ outside $\Omega$. Without loss of generality, we assume that $\|u\|_{N^{1,\Phi}(\Omega,V)}\leq1$. For any $x\in \bx$ and $k\in \nn$, we define 
	$$T_k u(x)=\vint_{B(x, 2^{-k} R)\cap\Omega} u\, d\mu.$$ 
	
	We first show that the limits
	$$\widetilde Tu=\lim_{k\rightarrow \infty} T_k u$$
	exist $\mH$-almost everywhere on $\bx$.
	It suffices to show that the function 
	$$\bar u=\sum_{k\geq 0} |T_{k+1}u -T_k u|+|T_0 u|$$
	belongs to $L^\Phi(\bx)$, since $\bar u\in L^\Phi(\bx)$ implies that $\bar u(x)<\infty$ for $\mH$-almost everywhere $x\in \bx$.
	Then it suffices to show that
	$$\|\bar u\|_{L^\Phi(\bx,V)}\leq \|T_0 u\|_{L^\Phi(\bx, V)}+\sum_{k\geq 0}\|T_{k+1}u- T_k u\|_{L^\Phi(\bx, V)}<\infty.$$
	
	Notice that $T_0 u(x)=\vint_{B(x, R)\cap\Omega} u\, d\mu=\vint_{\Omega} u\, d\mu$ for any $x\in \bx$, since $\Omega\subset B(x,R)$ for any $x\in\bx$. It follows from the Jensen inequality, upper codimension relation \eqref{upper} and $\Phi\in\Delta'$ that $$\|T_0 u\|_{L^{\Phi}(\bx, V)}\leq ac_1\|u\|_{L^{\Phi}(\Omega, V)},$$ where $a=\frac{1}{\Phi^{-1}(R^\theta)}$ and constant $c_1$ does not depend on $u$. Indeed, we have 
	\begin{align}
	\int_\bx\Phi\left(\frac{T_0 u}{ac_1\|u\|_{L^{\Phi}(\Omega, V)}}\right)\,d\mH &\leq \int_\bx\Phi\left(T_0\left(\frac{ u}{ac_1\|u\|_{L^{\Phi}(\Omega, V)}}
	\right)\right)\,d\mH\notag\\
	&\lesssim \int_\bx \frac{R^{-\theta}}{\mH(\bx)}\int_\Omega\Phi\left(\frac{u}{ac_1\|u\|_{L^{\Phi}(\Omega, V)}}\right)\,d\mu\,d\mH\notag\\
	&\lesssim R^{-\theta}\Phi\left(\frac{1}{c_1}\right)\Phi\left(\frac{1}{a}\right)\int_\Omega\Phi\left(\frac{u}{\|u\|_{L^{\Phi}(\Omega, V)}}\right)\,d\mu\leq 1,\notag
    \end{align}
    where the lase inequality follows by choosing the constant $c_1$ large enough.
    
	For any $k\geq 0$, it follows from the doubling property of $\mu$, the local $\Phi$-Poincar\'e inequality, $\Phi\in\Delta'\cap\nabla_2$ and the upper codimension relation \eqref{upper} that
	\begin{equation}\label{eq:long estimate}
		\|T_{k+1}u- T_k u\|_{L^\Phi(\bx, V)}
		\leq a_k\cdot c_2\|g_u\|_{L^\Phi(\Omega)},
	\end{equation}
	where $a_k=\frac{2^{-k}R}{\Phi^{-1}\left((2^{-k} R)^{\theta}\right)}$ and constant $c_2$ does not depend on $u$. Indeed, we have
	\begin{align}
		&\int_\bx\Phi\left(\frac{T_{k+1}u(x)-T_k u(x)}{a_k\cdot c_2\|g_u\|_{L^\Phi(\Omega)}}\right)\,d\mH(x)\notag\\
		&\leq \int_\bx\Phi\left(\frac{1}{a_k\cdot c_2\|g_u\|_{L^\Phi(\Omega)}} \cdot \vint_{B(x,2^{-k-1}R)\cap \Omega}|u(y)-u_{B(x,2^{-k}R)\cap \Omega}|\,d\mu(y)\right)\,d\mH(x)\notag\\
		&\lesssim \int_\bx\Phi\left(\frac{1}{a_k\cdot c_2} \cdot \vint_{B(x,2^{-k}R)\cap \Omega}|\frac{u(y)}{\|g_u\|_{L^\Phi(\Omega)}}-\frac{u_{B(x,2^{-k}R)\cap \Omega}}{\|g_u\|_{L^\Phi(\Omega)}}|\,d\mu(y)\right)\,d\mH(x)\notag\\
		&\lesssim \int_\bx\Phi\left(\frac{2^{-k}R}{a_k\cdot c_2}	\cdot\Phi^{-1}\left(\vint_{B(x,2^{-k}\lambda R)\cap\Omega}\Phi\left(\frac{g_u(y)}{\|g_u\|_{L^\Phi(\Omega)}}\right)\,d\mu(y)\right)\right)\,d\mH(x)\notag\\
		&\lesssim \int_\bx\Phi\left(\frac{2^{-k}R}{a_k\cdot c_2}\right)\vint_{B(x,2^{-k}\lambda R)\cap\Omega}\Phi\left(\frac{g_u(y)}{\|g_u\|_{L^\Phi(\Omega)}}\right)\,d\mu(y)\,d\mH(x)\notag\\
		&\lesssim
		\int_\bx\Phi\left(\frac{2^{-k}R}{a_k\cdot c_2}\right) \frac{(2^{-k}\lambda R)^{-\theta}}{\mH(B(x,2^{-k}\lambda R)\cap\bx)}\int_{B(x,2^{-k}\lambda R)\cap\Omega}\Phi\left(\frac{g_u(y)}{\|g_u\|_{L^\Phi(\Omega)}}\right)\,d\mu(y)\,d\mH(x)\notag\\
		&\lesssim 
		\int_{\Omega(2^{-k}\lambda R)}\Phi\left(\frac{g_u(y)}{\|g_u\|_{L^\Phi(\Omega)}}\right)\int_{B(y,2^{-k}\lambda R)\cap\bx} \Phi\left(\frac{2^{-k}R}{a_k\cdot c_2}\right) \frac{(2^{-k}\lambda R)^{-\theta}}{\mH(B(x,2^{-k}\lambda R)\cap\bx)}\,d\mH(x)\,d\mu(y)\notag\\
		&\lesssim \Phi\left(\frac{1}{c_2}\right)\Phi\left(\frac{2^{-k}R}{a_k}\right)(2^{-k}R)^{-\theta}\int_{\Omega(2^{-k}\lambda R)}\Phi\left(\frac{g_u(y)}{\|g_u\|_{L^\Phi(\Omega)}}\right)\,d\mu(y)\leq 1,\notag		
	\end{align}
	where $\Omega(r):=\{x\in \Omega: d(x, \bx)<r\}$ and the constant $c_2$ is large enough.
	
	Since $\sum_{k\geq 0}\frac{2^{-k}}{\Phi^{-1}\left(2^{-k\theta}\right)}<\infty$, combing the estimates of $\|T_0 u\|_{L^{\Phi}(\bx, V)}$ and $\|T_{k+1}u- T_k u\|_{L^\Phi(\bx, V)}$, we obtain that
	\begin{align*}
	\|\bar u\|_{L^\Phi(\bx,V)}&\lesssim \frac{1}{\Phi^{-1}(R^\theta)}\|u\|_{L^{\Phi}(\Omega, V)}+\sum_{k\geq 0}\frac{2^{-k}R}{\Phi^{-1}\left((2^{-k} R)^{\theta}\right)}\|g_u\|_{L^\Phi(\Omega)}\\
	&\lesssim \|u\|_{N^{1, \Phi}(\Omega, V)}<\infty.
    \end{align*}
	Thus, $\widetilde T u$ exists $\mH$-almost everywhere on $\bx$. Moreover, since $|\widetilde Tu|\leq \bar u$, we have 
	\[\|\widetilde Tu\|_{L^\Phi(\bx, V)}\leq \|\bar u\|_{L^\Phi(\bx,V)} \lesssim \|u\|_{N^{1, \Phi}(\Omega, V)}.\]
	
	The proof will be complete once we show  $\widetilde Tu=Tu$ on $\bx$. For this, it suffices to show that the eqaution \eqref{trace-defn} holds with $\widetilde Tu(x)$ for $\mH$-almost every $x\in \bx$. Set 
	$$E=\{x\in \bx: M_{\theta, \Phi} g_u(x)=\infty\ \text{and}\ T_k u(x)\rightarrow \widetilde Tu(x)\ \text{as}\ k\rightarrow \infty\}.$$ 
	Then Lemma \ref{fractional max operator} implies that $\mH(E)=0$. 
	
	For any  $0<r\leq R$, let $k_r\in \nn$ such that $2^{-k_r-1}R<r\leq 2^{-k_r}R$. Then it follows from the doubling property of $\mu$ and the local $\Phi$-Poincar\'e inequality
	that for any $x\in \bx\setminus E$ and $0<r\leq R$, 
	\begin{align*}
		\vint_{B(x, r)\cap\Omega} |u-\widetilde Tu(x)|\, d\mu
		&\leq \vint_{B(x, r)\cap\Omega} |u-T_{k_r}(x)|\, d\mu+|T_{k_r}(x)-\widetilde Tu(x)| \\
		&\lesssim \vint_{B(x, 2^{-k_r}R)\cap \Omega}|u-u_{B(x, 2^{-k_r}R)\cap \Omega}|\, d\mu +|T_{k_r}(x)-\widetilde Tu(x)| \\
		&\lesssim 2^{-k_r}R \cdot\Phi^{-1}\left(\vint_{B(x, 2^{-k_r}\lambda R)\cap \Omega} \Phi(g_u(x))\, d\mu(x)\right) +|T_{k_r}(x)-\widetilde Tu(x)|\\
		&\lesssim 2^{-k_r}R\cdot\Phi\left((2^{-k_r}R)^\theta\right) M_{\theta, \Phi} g_u(x)+ |T_{k_r}(x)-\widetilde Tu(x)|,
	\end{align*}
   where in the last inequality, we used the fact that $\Phi^{-1}\in \nabla'$ and $\Phi\in \Delta'$. Since $x\in \bx\setminus E$ and $k_r\rightarrow \infty$ as $r\to 0$, we have
	$$\vint_{B(x, r)\cap\Omega} |u-\widetilde Tu(x)|\, d\mu\to 0,\ \ \text{as}\ \ r\to 0.$$
	Hence \eqref{trace-defn} holds with $\widetilde Tu(x)$ for $\mH$-almost every $x\in \bx$. The proof is complete.
\end{proof}

\begin{remark}
	The condition $\sum_{k\geq 0}\frac{2^{-k}}{\Phi^{-1}\left(2^{-k\theta}\right)}<\infty$ is reasonable. Indeed, since the N-function $\Phi\in\Delta'\cap\nabla_2$, then $\Phi^{-1}(x)\geq C |x|^\frac{1}{p}$ for some $p>1$. Hence, if we choose $\theta<p$, then 
	$$\sum_{k\geq 0}\frac{2^{-k}}{\Phi^{-1}\left(2^{-k\theta}\right)}\leq C\sum_{k\geq 0}2^{-k(1-\frac{\theta}{p})}<\infty.$$
\end{remark}

\begin{proof}[Proof of Theorem \ref{thm:trace bdd-1}]
	This is a direct consequence of Theorem \ref{thm:trace bdd}.
\end{proof}

As a consequence of the proof of Theorem \ref{thm:trace bdd}, we obtain the following convergence result for traces of metric valued Orlicz-Sobolev space, which in particular gives Theorem \ref{thm:trace convergence1}. 

\begin{theorem}\label{thm:trace convergence}
	Suppose $\Omega\subset X$ is weakly $(\Phi,\theta)$-admissible for some N-function $\Phi\in\Delta'\cap\nabla_2$. Let $\{u_i\}\subset N^{1,\Phi}(\Omega,Y)$ be a sequence with uniformly bounded energy, that is, 
	$$\sup_{i\in \mathbb{N}}E^{\Phi}(u_i)<\infty.$$
	If $u_i$ converges to some $u\in N^{1,\Phi}(\Omega,Y)$ in $L^\Phi(\Omega,Y)$, then $Tu_i\to Tu$ in $L^\Phi(\bx,Y)$. Furthermore, two maps $u, v\in N^{1, \Phi}(\Omega, Y)$ have the same trace if and only if $d(u,v)\in N^{1, \Phi}(\Omega, \real)$ and has zero trace.
\end{theorem}

\begin{proof}
	For both assertions, embedding $Y$ isometrically into some Banach space $V$ if necessary, we may assume $Y=V$ is a Banach space.
	
	For the first claim, recall that in the proof of Theorem \ref{thm:trace bdd}, we proved that $Tf= \widetilde Tf$ for any $f\in N^{1, \Phi}(\Omega, V)$, where
	$$\widetilde Tf =\lim_{k\rightarrow \infty}T_k f.$$
	It follows from the estimate \eqref{eq:long estimate} that 
	\begin{align*}
		\|Tf-T_k f\|_{L^\Phi(\bx, V)} &\leq \sum_{j\geq k}\|T_{j+1} f-T_j f\|_{L^\Phi(\bx, V)}\\
		& \lesssim \sum_{j\geq k}\frac{2^{-j}R}{\Phi^{-1}\left((2^{-j} R)^{\theta}\right)}\|g_f\|_{L^\Phi(\Omega)}\\
		&\lesssim \frac{2^{-k}R}{\Phi^{-1}\left((2^{-k} R)^{\theta}\right)}\|g_f\|_{L^\Phi(\Omega)},
	\end{align*}
	where $g_f$ is the mimimal $\Phi$-weak upper gradient of $f$.
	
	Hence for any two maps $f, h\in N^{1,\Phi}(\Omega,V)$ and any $k\in \nn$, we have 
	\begin{align}\label{eq:difference of tf and th}
		\|Tf -Th\|_{L^\Phi(\bx, V)} 
		&\leq \|Tf-T_k f\|_{L^\Phi(\bx, V)}+\|Th-T_k h\|_{L^\Phi(\bx, V)}+\|T_kf-T_k h\|_{L^\Phi(\bx, V)} \notag\\
		&\lesssim \frac{2^{-k}R}{\Phi^{-1}\left((2^{-k} R)^{\theta}\right)}\left(\|g_f\|_{L^\Phi(\Omega)}+\|g_h\|_{L^\Phi(\Omega)}\right)+\|T_kf-T_k h\|_{L^\Phi(\bx, V)}, 
	\end{align}
	where $g_f$ and $g_h$ are minimal $\Phi$-weak upper gradients of $f$ and $h$, respectively.
	Notice that for any $x\in \bx$, we have
	$$T_k f(x)=\vint_{B(x, 2^{-k} R)\cap\Omega}f \, d\mu\ \ \text{and}\ \ T_k h(x)=\vint_{B(x, 2^{-k} R)\cap\Omega} h\, d\mu.$$
	Thus
	\begin{align*}
		\|T_kf-T_k h\|_{L^\Phi(\bx, V)}&\leq\|T_{k+1}h-T_k f\|_{L^\Phi(\bx, V)}+\|T_{k+1}h-T_k h\|_{L^\Phi(\bx, V)}\\
		&=: I_1+I_2.
	\end{align*}
	Using similar arguments as that in \eqref{eq:long estimate}, we obtain 
	$$I_2\lesssim \frac{2^{-k}R}{\Phi^{-1}\left((2^{-k} R)^{\theta}\right)}\|g_h\|_{L^\Phi(\Omega)}.$$
	For the estimate of $I_1$, it follows from the doubling property of $\mu$, the Jensen inequality and upper codimension relation \eqref{upper} that
	$$I_1\leq b_k\cdot c_3\|h-f\|_{L^\Phi(\Omega,V)},$$
	where $b_k=\frac{1}{\Phi^{-1}\left((2^{-k} R)^{\theta}\right)}$ and constant $c_3$ does not depend on $u$. Indeed, we have
	\begin{align}\label{eq:constant 2}
		&\int_\bx\Phi\left(\frac{T_{k+1}h(x)-T_kf(x)}{b_k\cdot c_3\|h-f\|_{L^\Phi(\Omega,V)}}\right)\,d\mH(x)\notag\\
		&\leq \int_\bx\Phi\left[\frac{1}{b_k\cdot c_3\|h-f\|_{L^\Phi(\Omega,V)}}\left(\vint_{B(x, 2^{-k-1} R)\cap\Omega}h(y)\,d\mu(y)-\vint_{B(x, 2^{-k} R)\cap\Omega}f(y)\,d\mu(y)\right)\right]\,d\mH(x)\notag\\
		&\lesssim \int_\bx\Phi\left(\frac{1}{b_k\cdot c_3}\right)\Phi\left(\vint_{B(x, 2^{-k} R)\cap\Omega}\frac{h(y)-f(y)}{\|h-f\|_{L^\Phi(\Omega,V)}}\,d\mu(y)\right)\,d\mH(x)\notag\\
		&\leq \int_\bx\Phi\left(\frac{1}{b_k\cdot c_3}\right)\vint_{B(x, 2^{-k} R)\cap\Omega}\Phi\left(\frac{h(y)-f(y)}{\|h-f\|_{L^\Phi(\Omega,V)}}\right)\,d\mu(y)\,d\mH(x)\notag\\
		&\lesssim \int_\bx\Phi\left(\frac{1}{b_k\cdot c_3}\right)\frac{(2^{-k}R)^\theta}{\mH(B(x, 2^{-k} R)\cap\bx)}\int_{B(x, 2^{-k} R)\cap\Omega}\Phi\left(\frac{h(y)-f(y)}{\|h-f\|_{L^\Phi(\Omega,V)}}\right)\,d\mu(y)\,d\mH(x)\notag\\
		&\lesssim \int_{\Omega(2^{-k} R)} \Phi\left(\frac{h(y)-f(y)}{\|h-f\|_{L^\Phi(\Omega,V)}}\right) \int_{B(y, 2^{-k} R)\cap\bx}\Phi\left(\frac{1}{b_k\cdot c_3}\right)\frac{(2^{-k}R)^\theta}{\mH(B(x, 2^{-k} R)\cap\bx)}\,d\mH(x)\,d\mu(y)\notag\\
		&\lesssim \Phi\left(\frac{1}{c_3}\right)\Phi\left(\frac{1}{b_k}\right)(2^{-k}R)^\theta \int_{\Omega(2^{-k} R)} \Phi\left(\frac{h(y)-f(y)}{\|h-f\|_{L^\Phi(\Omega,V)}}\right)\,d\mu(y)\leq 1,
	\end{align}
where $\Omega(r):=\{x\in \Omega: d(x, \bx)<r\}$ and the constant $c_3$ is large enough.

	Thus, the estimate \eqref{eq:difference of tf and th} can be rewritten as
	\begin{align}\label{eq:difference of tf and th result}
		\|Tf-Th\|_{L^\Phi(\bx, V)}
		&\lesssim  \frac{2^{-k}R}{\Phi^{-1}\left((2^{-k} R)^{\theta}\right)}\left(\|g_f\|_{L^\Phi(\Omega)}+\|g_h\|_{L^\Phi(\Omega)}\right)\notag\\
		&\quad\ +\frac{1}{\Phi^{-1}\left((2^{-k} R)^{\theta}\right)}\|f-h\|_{L^\Phi(\Omega,V)}.
	\end{align}
	The above inequality shows that if the sequence $u_i$ converges to $u$ in $L^\Phi(\Omega, V)$ and if the sequence has uniformly bounded energy, then $Tu_i$ converges to $Tu$ in $L^\Phi(\bx, V)$. Indeed, if we choose $f=u$ and $h=u_i$ in the above inequality, we know from the lower semicontinuity of energy (see \cite[Corollary 3.6.7]{HH19}, \cite[Theorem 2.4.1]{hkst12}) that the energy of $u$ is also bounded and hence the first term on the right-hand side of \eqref{eq:difference of tf and th result} can be made arbitrary small by choosing $k$ big enough. Once $k$ is fixed, the second term can be made small by choosing $i$ large.
	
	We now turn to the second claim and assume that $u, v\in N^{1, \Phi}(\Omega, V)$ have the same trace, i.e., $Tu(x)=Tv(x)$ for $\mH$-a.e. $x\in \bx$. We first show that $d(u,v)=|u-v|\in N^{1, \Phi}(\Omega)$. Since $|u-v|\leq |u|+|v|$, $|u-v|\in L^\Phi(\Omega)$. The minimal $\Phi$-weak upper gradient $g_{|u-v|}$ of $|u-v|$ is controlled by $g_u+g_v$, where $g_u$ and $g_v$ are minimal $\Phi$-weak upper gradients of $u$ and $v$. Indeed, for any rectifiable curve $\gamma$ connecting $x, y\in \Omega$, by triangle inequality, we have that
	\begin{align*}
		\big||u(x)-v(x)|-|u(y)-v(y)|\big|\leq |u(x)-u(y)|+|v(x)-v(y)| \leq \int_\gamma g_u+g_v\, ds.
	\end{align*}
	Thus, $|u-v|\in N^{1, \Phi}(\Omega)$. Since $Tu(x)=Tv(x)$ for $\mH$-a.e. $x\in \bx$, it follows from the definition of trace that for $\mH$-a.e. $x\in \bx$, we have
	\begin{align*}
		\lim_{r\rightarrow 0^+} \vint_{B(x, r)\cap\Omega} |u-v|\, d\mu&\leq \lim_{r\rightarrow 0^+} \vint_{B(x, r)\cap\Omega} |u-Tu(x)|\,d\mu+ |Tu(x)-Tv(x)|\\
		&\quad\quad\quad +\lim_{r\rightarrow 0^+} \vint_{B(x, r)\cap\Omega} |Tv(x)-v|\,d\mu=0.
	\end{align*}
	Hence $|u-v|$ has trace zero.
	
	For the converse, assume that $|u-v|\in N^{1, \Phi}(\Omega)$ has trace zero. Notice that for any $x\in \bx$ and any $y\in {B(x, r)\cap\Omega}$, we have that
	$$|Tu(x)-Tv(x)| \leq |Tu(x)-u(y)|+|u(y)-v(y)|+|v(y)-Tv(x)|.$$
	It follows from the definition of trace that for $\mH$-a.e. $x\in \bx$, we have
	\begin{align*}
		|Tu(x)-Tv(x)|&\leq \vint_ {B(x, r)\cap\Omega} |Tu(x)-u|\, d\mu+ \vint_ {B(x, r)\cap\Omega} |u-v|\, d\mu\\
		&\quad\quad\quad+ \vint_ {B(x, r)\cap\Omega} |v-Tv(x)|\, d\mu\rightarrow 0, \ \text{as}\ r\rightarrow 0.
	\end{align*}
	Thus, $u$ and $v$ have the same trace.
\end{proof}

\section{Solution to the Dirichlet problem}\label{sec:dirichlet problem}

In this section, we provide the proof of Theorem \ref{thm:main existence}, which is very similar to \cite[Theorem 1.4]{GW20} and \cite[Theorem 1.4]{GHWX21}. In the first step, we prove the following result on ultra-limits of subsequences of Orlicz-Sobolev maps, which extends \cite[Theorem 1.6]{GW20} and \cite[Theorem 4.3]{GHWX21}.

\begin{theorem}\label{thm:ultra-completion}
	Suppose $\Omega\subset X$ is a weakly $(\Phi,\theta)$-admissible domain for some N-function $\Phi\in\Delta_2\cap\nabla_2$ and $Y_\omega$ is an ultra-completion of the complete metric space $Y$. If $\{u_k\}\subset N^{1,\Phi}(\Omega,Y)$ is a bounded sequence, then, after possibly passing to a subsequence, the map $\phi(z):=[(u_m(z))]$ belongs to $N^{1,\Phi}(\Omega,Y_\omega)$ and satisfies
	$$E^\Phi(\phi)\leq \liminf_{k\to \infty}E^\Phi(u_k).$$
	Moreover, if $Tu_k$ converges to some map $\rho\in L^\Phi(\bx,Y)$ $\mH$-almost everywhere on $\bx$, then $T\phi=\iota\circ \rho$. 
\end{theorem}

\begin{proof}
	The proof is essentially contained in \cite[Proof of Theorem 1.6]{GW20} and we present it again for the convenience of the readers. After possibly passing to a subsequence, we may assume that
	$$E^\Phi(u_k)\to\liminf_{m\to\infty}E^\Phi(u_m)$$
	as $k\to \infty$.
	
	Fix $y_0\in Y$ and apply the Rellich-Kondrachov compactness Theorem \ref{thm:main theorem}. After possibly passing to a subsequence, there exist a complete metric space $Z=(Z,d_Z)$, a compact subset $K\subset Z$, and isometric embedding $\varphi\colon Y\to Z$ and $v\in N^{1,\Phi}(\Omega,Z)=M^{1,\Phi}(\Omega,Z)$ such that $\varphi(y_0)\subset K$ and
	$v_k:=\varphi\circ u_k$ converges in $L^\Phi(\Omega,Z)$ to $v$ as $k\to \infty$. 
	After passing to a further subsequence, we may assume that $v_k$ converges almost everywhere to $v$ on $\Omega$. Let $N\subset \Omega$ be a set of $\mu$-measure zero such that $v_k(z)\to v(z)$ for all $z\in \Omega\backslash N$.

	Define a subset of $Z$ by $B:=\{v(z):z\in \Omega\backslash N \}$.
	The map $\psi\colon B\to Y_\omega$, given by $\psi(v(z))=[(u_k(z))]$ when $z\in \Omega\backslash N$
	is well-defined and isometric by~\cite[Lemma 2.2]{GW20}. Since $Y_\omega$ is complete, there exists a unique extension of $\psi$ to $\overline{B}$, which we denote again by $\psi$. After possibly redefining the map $v$ on $N$, we may assume that $v$ has image in $\overline{B}$ and hence $v$ is an element of $N^{1,\Phi}(\Omega,\overline{B})$. Now, we define a mapping by $$\phi(z):=\psi(v(z))=[(u_k(z))]$$ 
	and then $\phi$ belongs to $N^{1,\Phi}(\Omega,Y_\omega)$ and by the lower semicontinuity of upper gradient energy (see \cite[Corollary 3.6.7]{HH19}, \cite[Theorem 2.4.1]{hkst12}) it satisfies
	\begin{align}\label{eq:5}
		E^\Phi(\phi)\leq E^\Phi(v)\leq \liminf_{k\to \infty}E^\Phi(v_k)=\liminf_{k\to \infty} E^\Phi(u_k).
	\end{align}
	
	It remains to prove the trace equality. Suppose $Tu_k$ converges to some map $\rho\in L^\Phi(\bx,Y)$ almost everywhere on $\bx$. Arguing as in \cite[Page 104]{GW20}, we can find compact subsets $C_1\subset C_2\subset \cdots \subset Y$, isometric embedding $\varphi\colon Y\to Z$ and $v\in N^{1,\Phi}(\Omega,Z)$ such that $v_k:=\varphi\circ u_k$ converges in $L^\Phi(\Omega,Z)$ to $v$ as $k\to \infty$. Furthermore, if we set $C=\bigcup_{l=1}^\infty C_l$, then aftering passing to a further subsequence if necessary we may assume that $v_k$ converges to $v$ almost everywhere on $\Omega$. Let $N\subset \Omega$ be a set of $\mu$-measure zero such that $v_k(z)\to v(z)$ for all $z\in \Omega\backslash N$. 
	
	Define a subset of $Z$ by
	$$B:=\{v(z):z\in\Omega\backslash N\}\cup \varphi(C).$$
	The map $\psi\colon B\to Y_\omega$ given by 
	\begin{equation*}
		\begin{cases}
			\psi(v(z))=[(u_k(z))] & \text{ if } z\in \Omega\backslash N,\\
			\psi(\varphi(x))=\iota(x)=[(x)] & \text{ if }x\in C,
		\end{cases}
	\end{equation*}
	is well-defined and an isometric embedding by \cite[Lemma 2.2]{GW20}. Since $Y_\omega$ is complete, there exists a unique isometric extension of $\psi$ to $\overline{B}$, which we denote again by $\psi$. After possibly redefining the map $v$ on $N$, we may assume $v\in N^{1,\Phi}(\Omega,\overline{B})$. The map $\phi(z):=\psi(v(z))=[(u_k(z))]$ then belongs to $N^{1,\Phi}(\Omega,Y_\omega)$ and satisfies \eqref{eq:5}. Moreover, by Theorem \ref{thm:trace convergence}, we have that $Tv_k=\varphi\circ Tu_k$ converges to $\varphi\circ \rho$ almost everywhere on $\bx$ and a subsequence of $Tv_k$ converges to $Tv$ almost everywhere. It thus follows that $Tv=\varphi\circ \rho$ and hence
	$$T\phi=\psi\circ Tv=\psi\circ\varphi\circ \rho=\iota\circ \rho.$$ 
	The proof is complete.
	
\end{proof}

With Theorem \ref{thm:ultra-completion} at hand, the proof of Theorem \ref{thm:main existence} is immediate. 
\begin{proof}[Proof of Theorem \ref{thm:main existence}]
Let $\phi\in N^{1,\Phi}(\Omega,Y)$ and let $\{u_k\}\subset N^{1,\Phi}(\Omega,Y)$ be an energy minimizing sequence with $Tu_k=T\phi$ for each $k$. Then by the characterization of trace from Theorem \ref{thm:trace convergence}, $h_k(x)=d(u_k(x),\phi(x))\in N^{1,\Phi}_0(\Omega)$.  Since $\sup_{k}E^\Phi(h_k)<\infty$, it follows from the global $\Phi$-Poincar\'e inequality \eqref{eq:global poincare inequality} that $\sup_{k}\|h_k\|_{L^\Phi(\Omega)}<\infty$. Hence
	$$\sup_{k} \left[ \| d(y_0,u_k(x))\|_{L^\Phi(\Omega)}+E^\Phi(u_k)\right]<\infty.$$
Thus $\{u_k\}$ is a bounded sequence in $N^{1,\Phi}(\Omega,Y)$. Let $Y_\omega$ be an ultra-completion of $Y$ such that $Y$ admits a 1-Lipschitz retraction $P\colon Y_\omega\to Y$. After possibly passing to a subsequence, we may assume by Theorem \ref{thm:ultra-completion} that the map $v(z):=[(u_k(z))]$ belongs to $N^{1,\Phi}(\Omega,Y_\omega)$ and satisfies $Tv=\iota\circ T\phi$ and 
$$E^\Phi(v)\leq \lim_{k\to \infty}E^\Phi(u_k).$$
Since $P\colon Y_\omega\to Y$ is a 1-Lipschitz retraction, the map $u:=P\circ v$ belongs to $N^{1,\Phi}(\Omega,Y)$ and satisfies $Tu=T\phi$ and $E^\Phi(u)\leq \lim_{k\to\infty}E^\Phi(u_m)$. The proof is complete.	
\end{proof}

\subsection*{Acknowledgements}

The author would like to thank her supervisor Prof.~Chang-Yu Guo for posing this question to her and for many useful conservations. This work is her master thesis.

\def\cprime{$'$} \def\cprime{$'$} \def\cprime{$'$}

\end{document}